\documentclass[11pt,a4paper]{article}
\usepackage{amsfonts,amsmath,amsthm}
\usepackage{amstext}
\usepackage[pdftex]{graphicx, color} 

\newtheorem{theorem}{Theorem}[section]
\newtheorem{lemma}[theorem]{Lemma}

\newtheorem{corol}[theorem]{Corollary}

\newenvironment{remark}%
  {\par\medbreak\refstepcounter{theorem}%
    \noindent\textbf{Remark~\thetheorem. }}%
  {\par\medskip}

\newcommand{\vz}[1]{\ensuremath{\mathbb{#1}}}

\newcommand{\R}{{\vz R}}

\newcommand{\B}{{\cal B}}

\DeclareMathOperator{\supp}{supp}

\long\def\drop#1{}

\let\e\varepsilon
\let\epsilon\varepsilon

\def\XXint#1#2#3{{\setbox0=\hbox{$#1{#2#3}{\int}$}
     \vcenter{\hbox{$#2#3$}}\kern-.5\wd0}}

\usepackage{geometry}                
\usepackage{graphicx}
\usepackage{amssymb}
\usepackage{epstopdf}
\usepackage{latexsym}
\usepackage{amssymb}
\usepackage{amsfonts}
\usepackage{bm}
\usepackage{epsfig}
\usepackage{bbm}
\usepackage{amsmath}

\usepackage{psfrag}
\usepackage{color}
\usepackage{multicol}
\usepackage{subfig}
\usepackage{appendix}
\usepackage{relsize}

\title{Deblurring of One Dimensional Bar Codes via Total Variation Energy Minimisation}

\author{Rustum Choksi
\footnote{Department of Mathematics and Statistics, McGill University, Burnside Hall, 805 Sherbrooke Street West, Montreal, Quebec, H3A 2K6, Canada, rchoksi@math.mcgill.ca}
  \qquad \, Yves van Gennip 
\footnote{Department of Mathematics, University of California Los Angeles, 520 Portola Plaza, Math Sciences Building 6363, Los Angeles, California, 90095, USA, y.v.gennip@gmail.com}}

\begin{document}

\maketitle

\begin{abstract} 
Using total variation based energy minimisation we address the recovery of a blurred (convoluted) one dimensional (1D) bar code. 
We consider functionals defined over all possible bar codes with fidelity to a convoluted  signal of a bar code, and regularised by total variation. Our fidelity terms consist  of the  $L^2$ distance either directly to the measured signal or preceded by deconvolution. Key length scales and parameters are the $X$-dimension of the underlying bar code, the size of the supports of the convolution and deconvolution kernels, and  the fidelity parameter. 
For all functionals, we establish parameter regimes (sufficient conditions) wherein the underlying bar code is the unique minimiser. 
We also present some numerical experiments suggesting that these  sufficient conditions are not optimal and the energy methods are  quite robust for significant blurring. 

\end{abstract} 


\medskip




\section{Introduction and notation}\label{sec:intro}
\begin{figure} 
\centerline{	{\includegraphics[width=2.3in]{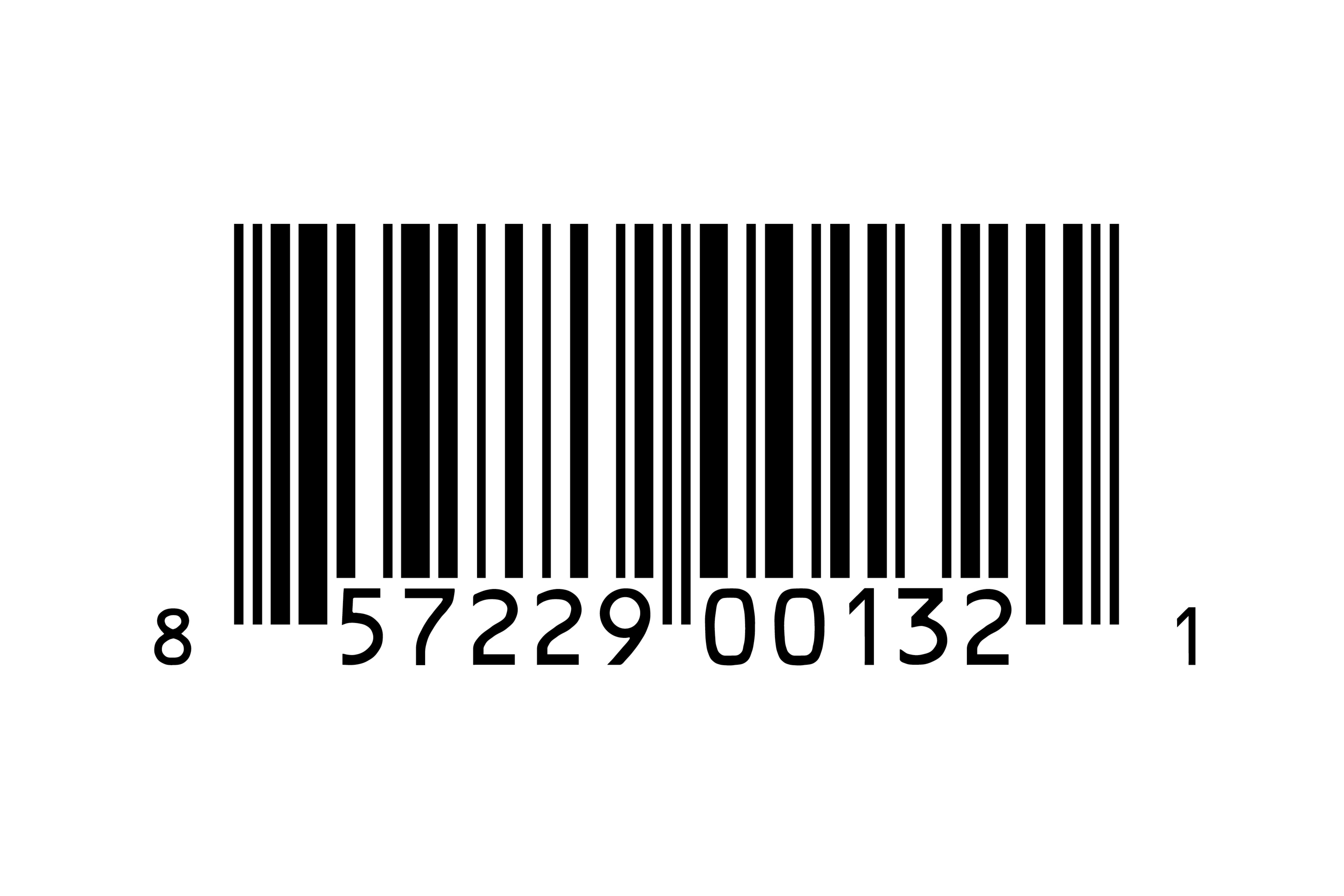}}\qquad {\includegraphics[width=2.3in]{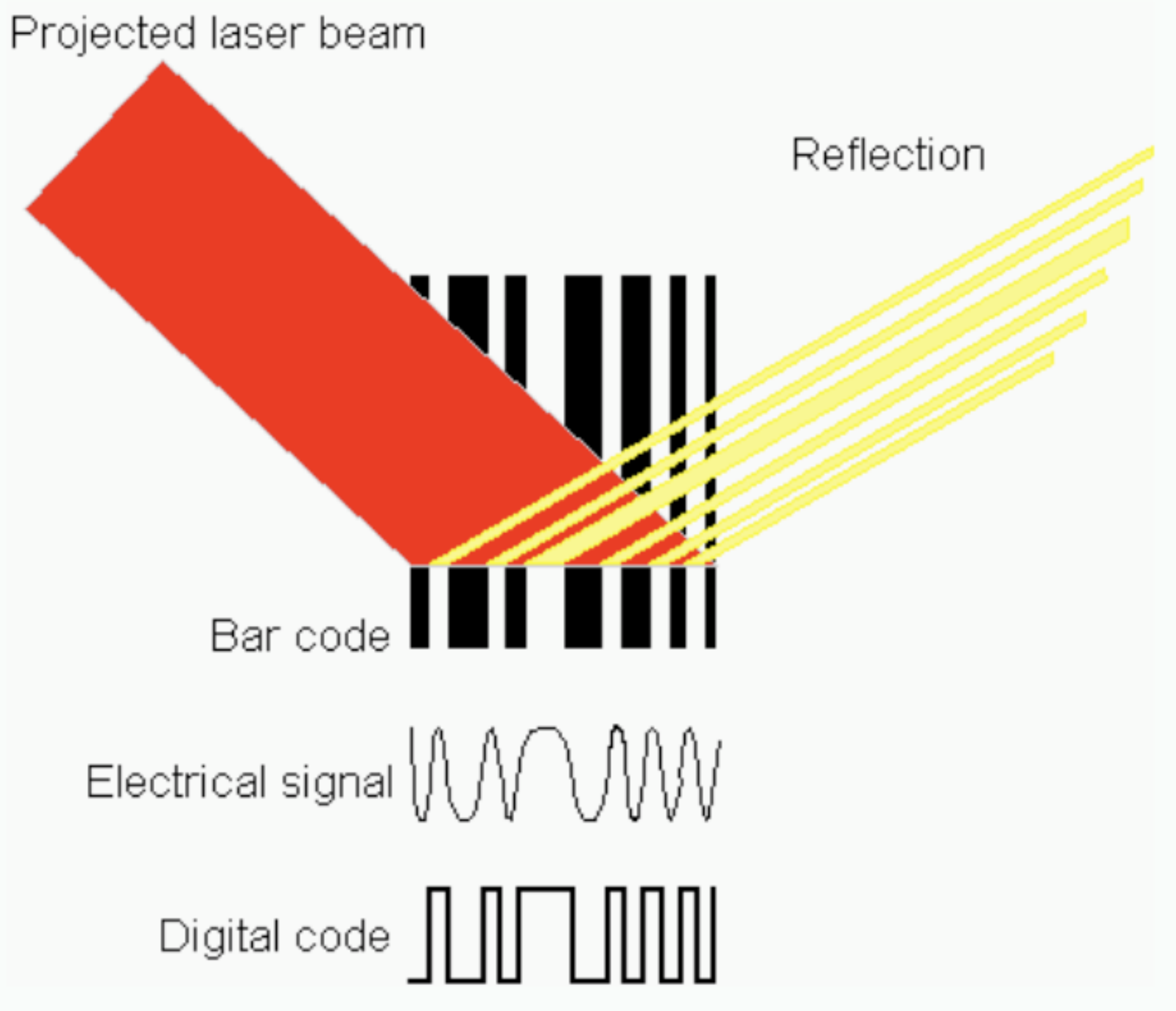}}}
\caption{Left: A standard 1D 12-digit Universal Product Code (UPC) bar code. The 12 numerical digits readable to the human eye are encoded in the bar spacings. \,  Right: A bar code scanner in which  the black bars absorb light, while the white bars reflect it. Photo diodes turn the reflected light into an electrical signal, which may be both blurred and noisy. This signal is then converted into  digital pulses. }
\label{Fig1}
\end{figure}

A one-dimensional bar code is a finite series of alternating black bars and white spaces with varying widths (Figure \ref{Fig1}).  The so-called 
{\it $X$-dimension} of the bar code is the width of the narrowest bar or space. 
A bar code scanner can use, for example, light detectors or photoreceptors similar to those used by a digital camera to 
pick up light and dark areas of the bar code, and produce a continuous signal associated with the darkness distribution 
across the span of the bar code ({c.f.} \cite{Palmer07}). 
One is left with an approximation  to the bar  code which depending on the scanner and the way 
in which the scan was taken (distance, vibrations etc.) can be blurred and noisy.
(see Figure  \ref{Fig1}). 
Thus in this article we 
deal with the general question: {\it Given a blurred  and possibly noisy {\it signal} $f$ associated with  a bar code, how can we deblur and denoise effectively to reconstruct the  bar code?}

Standard commercial deblurring techniques are often based upon edge detection, for example, finding local extrema of $f'$ which hopefully  correspond  to the discontinuities (i.e. interfaces)  of the original bar code. As noted in \cite{Esedoglu04}, this presents several difficulties: 
(i) the process is highly unstable to small changes in the signal, for example to the presence of noise;
(ii) points associated with local extrema of $f'$ are only crude approximations to the true locations of the bars with convolution tending to distort these points if the standard deviation of the convolution kernel ---or if the kernel is compactly supported the size of its support--- is comparable to the $X$-dimension of the bar code;
(iii) if the standard deviation or support size of the kernel is very large compared to the $X$-dimension, some edges in the bar code may have no corresponding extrema of $f'$ at all.

In this article, we take a different approach based upon energy minimisation involving a  total variation (TV) regularisation --  a method introduced by Rudin, Osher and Fatemi in \cite{RudinOsherFatemi92-2}.  One of the many advantages of this approach is  that energy minimisation  is a process {\it stable} with respect to small changes in the input signal.

In the context of bar code reconstruction, this type of energy minimisation was first proposed by Fadil Santosa, and analysed by Selim Esedo\=glu \cite{Esedoglu04} (see also \cite{Wittman04}). Here we take a similar approach but with an important difference: We are interested in directly testing the merits of energy minimisation for TV-based functionals by considering \emph{ans\"atze} for $f$ which involve convolutions of a bar code with certain blurring kernels. The functionals, blurring kernels, and admissible classes of bar codes possess certain length scale parameters. In \cite{Esedoglu04} Esedo\=glu shows existence of solutions for a variety of functionals and then proceeds to numerically test these functionals via  an algorithm which approximates the actual length scale parameter of the blurring kernel and uses this information to reconstruct the clean bar code signal. Here we fix these length scales as parameters and ask under what conditions  can we insure that energy minimisation gives back the underlying bar code.

We begin by introducing some notation. A bar code is given by a function $u$ 
of bounded variation ({\it c.f.} \cite{Giusti84}) with $\supp u \subset [0,1]$ taking on the values $0$ and $1$ a.e., i.e.  we consider a subset of the space $BV(\R; \{0, 1\})$. In particular the {\it general admissible set for bar codes} is 
\[
\B := \bigg\{ u \in BV({ \R}; \{0, 1\}) \,
: \,  u=0 \text{ a.e. in } [0,1]^c \bigg\}.
\]
Modulo a set of measure zero, the set $\{x\in[0,1]: u(x)=1\}$ consists of a finite number of disjoint non-empty intervals called bars.
We denote these non-empty bar intervals  by $b_i$ with length $|b_i|$. Similarly the intervals of nonempty  {\it white} spaces (i.e. intervals in $[0,1]$ where $u = 0$) are denoted by $w_i$. 
In this paper, signals will always be generated by bar codes whose $X$-dimension is bounded below. That is,  
we assume  there exists a constant $\omega >0$ such that 
the minimum width of these bars and spaces is {\it a priori} bounded below by  $\omega $.
Our space of {\it generating bar codes} is thus 
\[
\B_{\omega} := \bigg\{ u \in  \,\B 
:   \,\,\, \forall i \,\,\, \min\{ |b_i|, |w_i| \} \ge 
\omega \bigg\}.
\]

In some calculations it is useful, though harmless,  to assume that we know  whether the bar code starts (and/or ends) with a bar or a space. 
To this end,  we define for $i, j\in\{0,1\}$, the sets
\[
\B^{ij} := \Bigl\{ u\in \B: \text{there exists } x_1, x_2 \text{ such that } u=i \text{ on } [0,x_1] \text{ and } u=j \text{ on } [x_2, 1]\Bigr\},
\]
\[ \B_{\omega}^{ij} := \Bigl\{ u\in \B_{\omega}: \text{there exists } x_1, x_2 \text{ such that } u=i \text{ on } [0,x_1] \text{ and } u=j \text{ on } [x_2, 1]\Bigr\}. \]
As a useful addition to our terminology we will call the transition from a bar to a space or vice versa an \emph{interface}. This means that $\int_\R |u'|$ is equal to the number of interfaces in $u$ (where interfaces at $x=0$ and $x=1$ are included).
The notation $\chi_S$ will be used for the characteristic function of a given set $S$.

We approach the blurring via convolution on a length scale $\sigma >0$ (c.f. \cite{JosephPavlidis93, JosephPavlidis94, ShellhammerGorenPavlidis99, Esedoglu04}) with a symmetric unimodal kernel of unit mass. Different results hold for kernels of varying generality and we will discuss these kernels shortly. For the moment, let the kernel $\phi_\sigma$ denote a symmetric probability distribution  on $\R$ with {\it size}  (for example,  its standard deviation  or if compactly supported, half the size of its support) $\sigma$.  Given a bar code $z \in \B_{\omega}$,  we will assume observed signals of the form: 
\[ f_\sigma (x) \, := \, (\phi_\sigma \ast z) (x)  \, = \, \int_{-\infty}^\infty \phi_\sigma ( x - y) \, z (y) \, dy.\]
where $\sigma \ge 0$ is fixed (note that $\phi_0=\delta$, the Dirac delta distribution, and thus $f_0=z$).

For $u \in \B$:  
\begin{itemize}
\item We consider a fidelity term which compares $u$ directly with the observed signal:
\[
F_1(u) \, := \, \int_{\R} |u'| \, + \, \lambda \|u \, -\, f_\sigma \|_{L^2(\R)}^2.\]
\item ({\it deconvolution/deblurring functional})
Under the belief that the signal involves a convolution with a known kernel, we may incorporate this convolution into 
the structure of the fidelity term and consider: 
\[ F_2 (u) \, := \, \int_\R |u'| \, + \, \lambda \|\phi_\sigma * u \, - \,  f_\sigma\|_{L^2(\R)}^2. \]
If $\phi_\sigma$ is known, one may  ask as to the merits of energy minimisation as we could simply Fourier transform the observed signal  $f_\sigma$ and divide by the Fourier transform of $\phi_\sigma$ to recover the Fourier transform of $z$. However, this process is highly unstable with respect to small perturbations and in 
 practice, there is always some noise associated with the observed signal $ f_\sigma$. Directly solving back for $z$ is analogous to solving the heat equation backwards. Energy minimisation provides a stable numerical approach to  deblurring and indeed denoising.  


\item  ({\it blind deconvolution/deblurring functional})
Assuming no knowledge of $\sigma$, we may deconvolute with a kernel of similar type but with  
{\it size} $\rho$:
\[ F_3 (u)\, := \, 
\int_\R |u'| \, + \, \lambda \|\phi_\rho * u \,  - \,  f_\sigma\|_{L^2(\R)}^2.\]
\end{itemize} 
Note that $F_1$ and $F_2$ are special cases of 
$F_3$ with $\rho$ equal to $0$ and  $\sigma$ respectively,  and that in our notation,  the dependence of $F_i$ on $\lambda$, 
$z$, $\phi_\sigma$ and  $\phi_\rho$ is suppressed.
It is straightforward to see (cf. Lemma \ref{lemma0}), that for all parameters and all generating bar codes $z$, a minimiser of $F_i(u)$ over $u\in\B$  exists. 
In this article, we examine the following questions: 
{\it For what values of the parameters $\lambda, \omega, \sigma$ and $\rho$, is the minimiser of $F_i$ known, and in particular, when is it the underlying bar code $z$}? 
Our results consist of two parts:  First off we present a simple result that if $\lambda$ is sufficiently small, the unique minimiser is simply $u \equiv 0$. Part 2 deals with sufficient conditions for when the unique minimiser is the underlying bar code $z$. 

Note that if $\lambda$ is sufficiently small,  $u \equiv 0$ is the unique minimiser in $\B$ (i.e. an empty bar code). This is clearly the case if
\[
\lambda < \lambda_0 := 2 / \|f_\sigma \|_{L^2}^2,
\]
since any nontrivial bar code has a minimum total variation of $2$.
Sufficient and necessary conditions on $\lambda$ for $u=0$ to be the unique minimiser in $BV(\R)$ are given in the first parts of the Theorems \ref{thm:minimofF2},  \ref{thm:F4} and Corollary \ref{thm:F3} by exploiting a method adopted from \cite{Meyer01} wherein the following dual norm  associated with $BV$ is used:  
\begin{equation}\label{def:starnorm}  
\|f\|_\ast := \sup \left\{ \left|\int_\R f v\right| : v \in L^1(\R) \text{ and } \int_\R |v'| \leq 1\right\}.
\end{equation}
This threshold for a trivial minimiser is 
given by $\lambda \leq \frac1{2 \|\phi_\rho*f_\sigma\|_\ast}$, $\rho \geq 0$ (parts 1 of the theorems below) and is lower than $\lambda_0$. 
To see this note that for any  $f\in L^2(\R; [0, 1])$ with compact support we may take 
$I$ to be a bounded subset of $\R$ such that $\supp \phi_\rho*f \subset I$. Since $v=\frac12 \chi_I$ is an admissible function in the definition of $\|\phi_\rho*f\|_*$  we have 
\[
\|\phi_\rho*f\|_* \geq \frac12 \int_I \phi_\rho*f = \frac12 \int_\R f  \geq \frac12 \int_\R f^2.
\]
Therefore
\[
\frac1{2 \|\phi_\rho*f\|_*}\leq \left(\int_\R f^2\right)^{-1} <  \frac2{\|f\|_{L^2}^2}.
\]
For $\lambda$ between this threshold and $\lambda_0$, $u=0$ is the unique minimiser in $\B$ but not in $BV(\R)$.


The case where $z$ is the unique minimiser is more subtle and depends critically on both the size and particular 
nature of the blurring kernel. We make two assumptions here: 
\begin{enumerate}
\item \label{item:Xdimassump} We restrict our attention to kernels with compact and {\it small}  (with respect to the $X$-dimension $\omega$) support. 
\item We further consider unimodal, symmetric kernels in $\mathcal{K}$ defined below 
and  
explicit regimes are computed using a prototype of such a kernel, the hat function defined by 
  \[
\hat \phi_{\sigma}(x) := \left\{ \begin{array}{ll} (1\, - \, {x}/{\sigma})/{\sigma} & \text{if } 0\leq x<\sigma,\\  
(1\,+ \, {x}/{\sigma})/{\sigma} & \text{if } -\sigma < x \leq 0,\\ 0 & \text{if }|x| \geq \sigma. \end{array}\right.
\]
The class $\mathcal{K}$ is defined by 
\begin{align}\label{eq:Ksigma}
\mathcal{K} &:= \left\{ \phi_\sigma \in L^1(\R): \exists \sigma>0 \,\,\, \phi_\sigma (x) = p(-x, \sigma) \chi_{[-\sigma, 0]}(x) + p(x, \sigma) \chi_{[0, \sigma]}(x)\right. \notag\\&\hspace{0.7cm} \left. \text{ for a non-negative function } p: [0, \sigma]\times(0, \infty) \text{ monotonically decreasing}  \right.\notag \\ &\hspace{0.7cm} \left.  \text{ in } x, \text{ and } \int_0^\sigma p(x, \sigma)\, dx = \frac12\right\}. 
\end{align}
We will consistently use a subscript as in $\phi_\sigma$ to indicate the value of the parameter $\sigma$ in the definition of $\mathcal{K}$, i.e. $\phi_\rho \in \mathcal{K}$ is in the subset of $\mathcal{K}$ where $\sigma=\rho$.

Note in particular that if $\phi_\sigma\in L^1(\R)$, we have $\phi_\sigma*u\in L^2(\R)$ for all $u\in \mathcal{B}$.

It is also convenient to consider a subclass of  $\mathcal{K}$ 
\begin{align}
\mathcal{K}_3 &:= \left\{ \phi_\sigma \in \mathcal{K}: \phi_\sigma \in C_c(\R), p \text{ continuously differentiable in } \sigma  \text{ and } (\ref{eq:condpsit}) \text{ holds}\right\},\notag   
\end{align}
where the non-obvious condition (\ref{eq:condpsit})  is 
\begin{align}\label{eq:condpsit}
&\hspace{1cm} \forall \tau \in (0, \sigma], \,\, \forall c \geq 2 \sigma, \,\,  \forall x\in [0,c]: \notag\\
&\mathcal{J}(\sigma, \tau, x, c) := \int_0^\tau \int_{x-c}^x \frac{\partial}{\partial \tau} p(y, \tau) \left[ \phi_\sigma(y-w)+\phi_\sigma(y+w)\right]\, dw\, dy \leq 0.
\end{align}

This condition insures a certain monotonicity property (c.f. Lemma \ref{lem:decreasing}) of double convolutions with bar codes. As we show in Appendix~\ref{app:psitaufrhononpos} a sufficient condition for (\ref{eq:condpsit}) to be satisfied is if for each $\tau \in (0, \sigma]$
\begin{enumerate}
\item either $\displaystyle \frac{\partial}{\partial \tau} p(x, \tau)$ is monotonically increasing in $x$ and $\mathcal{J}(\sigma, \tau, 0, c) \leq 0$ for all $c\geq 2\sigma$,
\item or $\displaystyle \frac{\partial}{\partial \tau} p(x, \tau)$ is monotonically decreasing in $x$ and $\mathcal{J}(\sigma, \tau, \frac{c}2, c) \leq 0$ for all $c\geq 2\sigma$.
\end{enumerate}
In the same appendix we show that the hat function $\hat \phi_\sigma$ satisfies condition (a) above for each $\tau\in (0, \sigma]$. In practical applications using kernels that do not satisfy either (a) or (b), condition (\ref{eq:condpsit}) can be tested numerically.
%
\end{enumerate}
Assumption~\ref{item:Xdimassump} is rather important and indeed restrictive as it limits the possible effect of blurring. 
For $F_2$ and $F_3$ we require $\sigma, \rho \leq \omega/2$, insuring no interactions between neighbouring bars. For $F_1$, the condition is slightly less restrictive, namely  $\sigma \leq \omega$. 
The second assumption, particularly, the specification of the hat function is more for convenience. 
A crucial step in our direct and rather {\it brute-force} approach is to assume a minimiser with a certain structure and directly construct competitors which differ on a set bar or space.  
For this step, one can explicitly calculate a regime wherein such a competitor  exists, and for simplicity we have performed the calculations for the hat function (which were greatly simplified by the use of Maple). 
We discuss modifications for other kernels in Remark \ref{other} below. Let us now state our results.

\bigskip

\begin{theorem}\label{thm:minimofF2}
 The following hold: 
\begin{enumerate}
\item\label{item:F1trivminoverBV}  Let $\phi_\sigma \in \mathcal{K}$. $ u \equiv 0$ is the unique minimiser for $F_1$  over $BV(\R)$ iff $\|f_\sigma\|_\ast \leq \frac1{2\lambda}$.
\item\label{item:F1zminoverB} Let $\phi_\sigma = \hat \phi_\sigma$ and $z\in \B_{\omega}$. If $\sigma$ and $\lambda$ satisfy
\begin{equation}\label{eq:conditionsigmalambda}
\sigma \leq \omega \qquad {\rm and} \qquad \frac23 \sigma + \frac2{\lambda} < \omega,
\end{equation}
then $u = z$ is the unique minimiser of $F_1$ over $\B$.
\end{enumerate}
\end{theorem}


\begin{theorem}\label{thm:F4}
 The following hold: 
\begin{enumerate}
\item\label{item:F3trivminoverBV}  Let $\phi_\sigma \in \mathcal{K}$. $ u \equiv 0$ is the unique minimiser for $F_3$  over $BV(\R)$ iff $\|\phi_\rho*f_\sigma\|_\ast \leq \frac1{2\lambda}$. 
\item\label{item:F3zminoverB} Let $\phi_\sigma = \hat \phi_\sigma$ and $z
\in \B_{\omega}^{ij}$ for some $i, j\in \{0,1\}$. Let $\sigma \leq \rho \leq \frac{\omega}2$.
If $\lambda$, $\rho$, and $\sigma$ satisfy
\begin{equation}\label{eq:conditiononrhosigmalambda}
\frac2\lambda + \frac1{15 \rho^2} \Bigl(-\sigma^3 + 5 \rho \sigma^2 + 17 \rho^3\Bigr) < \omega,
\end{equation}
then $u=z$ is the unique minimiser of $F_3$ over $\B^{ij}$.
\end{enumerate}
\end{theorem}
Note that the left hand side of (\ref{eq:conditiononrhosigmalambda}) is increasing as a function of (real and positive) $\rho$ and $\sigma$.
By taking $\rho = \sigma$ in Theorem \ref{thm:F4}, we obtain 

\begin{corol}\label{thm:F3} The following hold: 
\begin{enumerate}
\item\label{item:F2trivminoverBV}  Let $\phi_\sigma \in \mathcal{K}$. $ u \equiv 0$ is the unique minimiser for $F_2$ over $BV(\R)$ iff \\
$\|\phi_\sigma*f_\sigma\|_* \leq \frac1{2\lambda}$. 
\item\label{item:F2zminoverB} Let $\phi_\sigma = \hat \phi_\sigma$ and $z
\in \B_{\omega}^{ij}$ for some $i, j\in \{0,1\}$. If  $\sigma \leq \frac{\omega}2$ and $\lambda>0$ satisfy   
\begin{equation}\label{eq:condsigmalambdaconvolutions}
\frac2\lambda + \frac{21}{15} \sigma < \omega, 
\end{equation}
then $u=z$ is the unique minimiser of $ F_2$ over $\B^{ij}$.
\end{enumerate}
\end{corol}

\begin{remark}\label{other} {\bf Extensions to other kernels}
We remark on extensions of parts 2 of Theorem~\ref{thm:minimofF2}, Theorem~\ref{thm:F4},  and Corollary~\ref{thm:F3} to more general kernels in $\mathcal{K}$. Their proofs consist of two steps: (i) First is to  establish that any minimiser of $F_1$ or  $F_3$ distinct from $z$ must have strictly less interfaces than $z$ (note that this is trivially satisfied for $F_2$ since $z$ uniquely minimises the fidelity term). 
For $F_1$, this is a simple consequence of vanishing first variation (Lemma \ref{lem:firstvar-F1}), and holds for any kernel in $\mathcal{K}$. For $F_3$, the consequences of vanishing first variation are more involved (c.f. Lemmas \ref{lem:levelhalfrhosigmaz} - \ref{lem:F3firstvar}), and an important ingredient is  
a monotonicity property of  double convolutions (c.f. Lemma \ref{lem:decreasing}) which is responsible for condition (\ref{eq:condpsit}). 
Thus, this step holds for any kernel in  $\mathcal{K}_3$.   
(ii) The second step involves  the explicit parameter regimes and is the reason why we have conveniently  adopted the hat function. 
Here we show that there cannot exist a minimiser for $F_1$, $F_2$, or $F_3$ 
 with fewer interfaces than $z$ if (\ref{eq:conditionsigmalambda}), (\ref{eq:condsigmalambdaconvolutions}), or (\ref{eq:conditiononrhosigmalambda}) holds respectively:  If a minimiser $u_0$ does have fewer interfaces than $z$, then there exists an interval (bounded below in length by $\omega$) on which $z$ has a bar and $u_0$ a space or vice versa. We then contradict our assumption by explicitly modifying $u_0$ on this interval to achieve  lower energy. This last step requires some straightforward but tedious calculations. Maple has been a great help in performing the many integrations necessary involving the hat function $\hat \phi_\sigma$. 
This step can be reproduced for any $\phi_\sigma \in \mathcal{K}$ with different parameter regimes for each specific choice of kernel $\phi_\sigma$ as a result. Specifically,   the calculations in the proof of Theorem~\ref{thm:minimofF2}, part~\ref{item:F1zminoverB} after (\ref{eq:twointegrals}) or the calculations in the proof of Lemma~\ref{lem:F3partresult} after (\ref{eq:F4inequality}) respectively need to be redone for the new kernel. 

Collecting the conditions necessary for steps (i) and (ii) we find that it is in principle possible to get results of the form of those in Theorem~\ref{thm:minimofF2}, part~\ref{item:F1zminoverB} and Corollary~\ref{thm:F3}, part~\ref{item:F2zminoverB} for $F_1$ and $F_2$ respectively for any $\phi_\sigma \in \mathcal{K}$. 
Similarly a result for $F_3$ as the one in Theorem~\ref{thm:F4}, part~\ref{item:F3zminoverB} can be obtained for any $\phi_\sigma\in \mathcal{K}_3$. 
(Note that Corollary~\ref{thm:F3}, part~\ref{item:F2zminoverB} can be obtained either as a consequence of Theorem~\ref{thm:F4}, part~\ref{item:F3zminoverB} or using Lemma~\ref{lem:F3partresult}. The latter option allows for less restrictions on $\phi_\sigma$ in Corollary~\ref{thm:F3} than in Theorem~\ref{thm:F4}.)
\end{remark}

We conclude this section with a few comments on the results,  their interpretations and limitations. 
Our {\it brute force}  arguments  are based upon explicit calculations and are as such limited to a blurring kernel with {\it small} (with respect to $\omega$) support.  
Indeed, this is  most discouraging  for the deconvolution functionals $F_2$ and $F_3$ where one 
 one would expect  the regime of acceptable $\sigma$ to extend far beyond the $X$-dimension of the underlying bar code. 
 Numerical experiments (see Section \ref{numerics}) support this conjecture. 
The conditions that $\sigma, \rho \leq w/2$ for $F_2$ and $F_3$ may seem particularly alarming since the analogous condition for $F_1$ is simply $\sigma \leq \omega$. However, note that the other condition (\ref{eq:conditiononrhosigmalambda}) also puts a restriction on the size of $\sigma$ which is essentially of the same form. 
For $F_2$, one could leave out the condition $\sigma \leq \frac\omega2$ and not change the principles of the proof, but  many more orderings in the computation of the integrals become possible (see Remark~\ref{rem:roleofconditions}) and many more calculations need to be done in the proof of Lemma~\ref{lem:F3partresult}. Since we still have condition (\ref{eq:condsigmalambdaconvolutions}) on $\sigma$ in place 
(these extra calculations can only replace (\ref{eq:condsigmalambdaconvolutions}) by a stricter condition on $\sigma$, if anything changes at all) it is doubtful that much can be won by leaving out the condition $\sigma \leq \frac\omega2$. For $F_3$ the conditions $\rho, \sigma \leq \frac\omega2$ are vital to our proofs via Lemmas~\ref{lem:levelhalfrhosigmaz},~\ref{lem:decreasing},~\ref{lem:F3firstvar}.

The numbers in condition (\ref{eq:conditiononrhosigmalambda})  may seem a little strange. They are simply a consequence of the direct calculations with the hat function. As we have remarked, these calculations can be repeated for other kernels in $\mathcal{K}$.  
Note that, taking $\rho = \omega/2$, the condition implies that the bar code is always recoverable 
for any $\sigma \leq \omega/2$, provided $\lambda > 20 / 3\omega$.
 
Another surprising condition might be $\sigma \leq \rho$ in Theorem~\ref{thm:F4}, part~\ref{item:F3zminoverB}. In general if $\rho <  \sigma$ we do not expect $z$ to be a minimiser of $F_3$ over $\B$ (or $\B^{ij}$), even if $\rho, \sigma \leq \frac\omega2$ and (\ref{eq:conditiononrhosigmalambda}) are satisfied. A counter example in this case is given by $z=\chi_{[0.425, 0.575]}$ and $u = \chi_{[0.425, 0.4999]} + \chi_{[0.5001, 0.575]}$ with $\rho=0.05$ and $\sigma=0.06$. The fidelity term in $F_3(u)$ is smaller than the fidelity term in $F_3(z)$ ($\|\phi_\rho*u-\phi_\sigma*z\|_{L^2(\R)}^2 \approx 2.378 \cdot 10^{-4}$ and $\|\phi_\rho*z-\phi_\sigma*z\|_{L^2(\R)}^2 \approx 2.407 \cdot 10^{-4}$) and thus for $\lambda$ large enough $z$ will not be the minimiser of $F_3$ in $\B$. In this particular case the difference is small and so in practical applications where $\lambda$ is not too large it might not cause problems, since then the energetic cost $2$ for two extra interfaces will be much higher than the gain in the fidelity term. 
Extra conditions on the parameters in the case $\rho  <  \sigma$ might ensure that $u=z$ is the minimiser for $F_3$. The above example suggests that an upper bound for $\lambda$ may be in order.
In fact, numerical simulations in  Section \ref{numerics} show that the regime $\rho < \sigma$ poses no problem for suitable {\it midrange}  choices of $\lambda$. Indeed,  they suggest that  fixing $\rho$ comparable with the $X$-dimension and minimising $F_3$ works well for $\sigma$ up to twice the $X$-dimension. 

We know however that in the degenerate case $0=\rho < \sigma$, 
Theorem~\ref{thm:minimofF2} assures that we have $u=z$ as minimiser if both conditions in (\ref{eq:conditionsigmalambda}) are satisfied, without an upper bound on $\lambda$. Why the second condition in (\ref{eq:conditionsigmalambda}) is the correct degenerate form of (\ref{eq:conditiononrhosigmalambda}) can be seen by recognizing their common source (\ref{eq:conditiononrhosigmalambda2}).





Finally, we note that there is a wealth of work on total variation energy minimisation for image analysis. While our results are for rather simple one dimensional images, we feel they are novel in that the 1D bar code setting entails an image deblurring  problem of contemporary interest  yielding  very precise results, and we are unaware of any general method  for analogous deblurring functionals which would yield similar results. 
In addition to geometric simplicity due to its binary nature (the simplest case of what is called \emph{Quantum TV} in \cite{ShenKang07}), the bar code problem is different from many other imaging problems in that there is a known {\it a priori} lower bound on 
the length scale of the structures in the image (via the $X$-dimension). An analytically deeper study entails deblurring of 
2D bar codes \cite{ChoksiGennipOberman10}.

\bigskip

\section{Proofs of the Theorems} \label{sec:proofs}

\subsection{Existence and the trivial minimiser}

\begin{lemma}\label{lemma0}
Let $\phi_\sigma \in \mathcal{K}$ and fix $z\in \B_\omega$ ($\B_\omega^{ij}$) and $\lambda, \sigma, \rho > 0$. Then minimisers for $ F_1, F_2$ and $ F_3$ over $\B$ ($\B^{ij}$) exist.
\end{lemma}
\begin{proof}
The proof is a simple application of the {\it direct method in the calculus of variations} and  follows along the same lines for all these functionals. For completeness,  we present it for  $F_1$ and $z\in \B_\omega$. Let $\{u_n\}$ be a minimising sequence for $F_1$ in $\B$, then we can assume every $u_n$ has bounded $L^1$-norm and bounded BV measure. Therefore, by \cite[\S 5.2.3 Theorem 4] {EvansGariepy92}, there exists $u\in BV([0,1])$ such that $u_n \to u$ in $L^1([0,1])$. Since the $u_n$ only take values $0$ and $1$ (and $0$ a.e. in $[0,1]^c$), so does $u$. 
 Thus  
$u \in \B$. 

The total variation is lower semicontinuous under $L^1$ convergence \cite[Theorem 1.9]  {Giusti84} and under the special conditions that the functions only take values $0$ and $1$, so is the $L^2$ norm, therefore we conclude via the direct method in the calculus of variations that $u$ is a minimiser for $F_1$. For the functionals $F_2$ and $F_3$ we use in addition, that the functional $u\mapsto \phi_\sigma * u$ is continuous under $L^1$ convergence, for any $\sigma>0$. If we replace $\B_\omega$ and $\B$ by $\B_\omega^{ij}$ and $\B^{ij}$ respectively the proof does not change.

\end{proof}
Note that we have not used the fact that $\phi_\sigma$ is symmetric and unimodal with compact support in the above. We only need continuity of $u\mapsto \phi_\sigma * u$ under $L^1$ convergence.

\bigskip

Next we recall a result about convolutions
whose proof follows directly from Fubini's Theorem.

\begin{lemma}\label{eq:intconv}
Let $f, g, h \in L^2(\R)$ such that $f(-x)=f(x)$, then
\[
\int_\R \big[(f*g)\cdot h\big] = \int_\R \big[g \cdot(f*h)\big].
\]
\end{lemma}

Parts 1 of Theorems \ref{thm:minimofF2}, \ref{thm:F4} and Corollary \ref{thm:F3} follow directly from the following lemma. 
\begin{lemma}\label{lem:trivialminimiser}
Let $i\in\{1, 2, 3\}$, $\phi_\sigma \in \mathcal{K}$, and $\lambda>0$,
then the following two statements are equivalent:
\begin{enumerate}
\item \label{item:trivialminimiser} $u=0$ is the unique minimiser of $F_i$ over $BV(\R)$. 
\item \label{item:starnorms} \begin{enumerate}
	\item If $i=1$, $\|f_\sigma\|_* \leq \frac1{2\lambda}$.
	\item If $i=2$, $\|\phi_\sigma*f_\sigma\|_* \leq \frac1{2\lambda}$.
	\item If $i=3$, $\|\phi_\rho*f_\sigma\|_* \leq \frac1{2\lambda}$.
	\end{enumerate}
\end{enumerate}
\end{lemma}
\begin{proof}
The idea of the proof is similar to that in \cite[\S1.14, Lemma 4]{Meyer01}. Note that for general $u\in BV(\R)$ we cannot conclude that $\phi_\sigma*u \in L^2(\R)$. For functions $u$ and parameters $\sigma$ or $\rho$ for which this fails we set $F_i(u)=\infty$, $i\in \{2, 3\}$.

First let $i=1$. We first prove \ref{item:trivialminimiser} $\Longrightarrow$ \ref{item:starnorms}. Assume $u=0$ is the unique minimiser of $F_1$ in $BV(\R)$.
This is equivalent to, for all $h\in BV(\R)$ with $h\neq 0$, 
\begin{equation}\label{eq:uisatrivialminimiser}
\lambda \|f_\sigma\|_{L^2(\R)}^2 < \int_\R |h'| + \lambda \|h-f_\sigma\|_{L^2(\R)}^2 = \int_\R |h'| + \lambda \|f_\sigma\|_{L^2(\R)}^2 + \lambda \|h\|_{L^2(\R)}^2 - 2 \lambda \int_\R f_\sigma h.
\end{equation}
Because this holds for all $h\in BV(\R)$, by rescaling $h$ we can rewrite this as
\begin{equation}\label{eq:uisatrivialminimiser2}
2 \lambda \e \int_\R f_\sigma h < |\e| \int_\R |h'| + \lambda \e^2 \|h\|_{L^2(\R)}^2,
\end{equation}
for all $\e\in\R$ and all $h\in BV(\R)$. Dividing by $\epsilon$, taking the limit $\e\to 0$, and recognizing that $\e$ can be positive and negative, we find that (\ref{eq:uisatrivialminimiser2}) implies
\begin{equation}\label{eq:equivstatement}
\left| \int_\R f_\sigma h \right| \leq \frac1{2\lambda} \int_\R |h'|, \qquad \text{for all } h\in BV(\R).
\end{equation}
Now per definition
\[
\|f_\sigma\|_* = \underset{v\in L^1(\R), \int_\R |v'| \leq 1}{\sup}\, \left|\int_\R f_\sigma v\right| 
\leq \frac1{2\lambda},
\]
where the inequality follows by taking the supremum in (\ref{eq:equivstatement}) over all\\
 $h\in \left\{v\in L^1(\R): \int_{\R} |v'| \leq 1\right\} \subset BV(\R)$.

To prove \ref{item:starnorms} $\Longrightarrow$ \ref{item:trivialminimiser} let $\|f_\sigma\|_* \leq \frac1{2\lambda}$. Then for all $v \in L^1(\R)$ satisfying $\int_\R |v'| \leq 1$ we have
\[
\left|\int_\R f_\sigma v \right| \leq \frac1{2\lambda},
\]
from which it follows that for all $h\in BV(\R)$
\[
\left| \int_\R f_\sigma \frac{h}{\int_\R |h'|}\right| \leq \frac1{2\lambda}.
\]
Inequality (\ref{eq:equivstatement}) now follows.

We proved above that (\ref{eq:uisatrivialminimiser}) implies (\ref{eq:equivstatement}). On the other hand we see that inequality (\ref{eq:equivstatement}) implies for $h\neq 0$
\[
\int_\R |h'| + \lambda \|f_\sigma\|_{L^2(\R)}^2 + \lambda \|h\|_{L^2(\R)}^2 - 2 \lambda \int_\R f_\sigma h \geq \lambda \|f_\sigma\|_{L^2(\R)}^2 + \lambda \|h\|_{L^2(\R)}^2 > \lambda \|f_\sigma\|_{L^2(\R)}^2,
\]
and thus inequality (\ref{eq:equivstatement}) is equivalent to (\ref{eq:uisatrivialminimiser}). 
This proves the result for $i=1$.

$F_2$ is a special case of $F_3$ (with $\rho=\sigma$). For $i=3$ we can derive a statement analogous to inequality (\ref{eq:equivstatement}), with $h$ on the left hand side replaced by $\phi_\rho*h$.  Having $u=0$ as unique minimiser of $F_3$ in $BV(\R)$ is equivalent to
\[
\left| \int_\R f_\sigma \cdot \phi_\rho*h \right| \leq \frac1{2\lambda} \int_\R |h'|, \qquad \text{for all } h\in BV(\R).
\]
By Lemma~\ref{eq:intconv} we recognise that
\[
\int_\R f_\sigma \cdot \phi_\rho*h = \int_\R \phi_\rho * f_\sigma \cdot h
\]
and the result follows as before.
\end{proof}

Although we assume $\phi_\sigma\in \mathcal{K}$ in the proof above because that is the most general class of kernels we consider, we only use the symmetry and integrability of $\phi_\sigma$.

\subsection{Proof of Theorem~\ref{thm:minimofF2}, part~\ref{item:F1zminoverB}}

We now turn our attention from the trivial minimiser $u=0$ to $u=z$ as minimiser. First,   we present  some elementary results.


\begin{lemma}\label{lem:onebarconvolution}
Let $\phi_\sigma \in \mathcal{K}$, $a<b$, $\sigma \leq b-a$ and $z=\chi_{[a, b]}$, then
\[
\left\{ x\in \R : f_\sigma(x) = \frac12\right\} = \{a, b\} \quad \text{and} \quad \left\{ x\in \R : f_\sigma(x) \geq \frac12\right\} = [a, b].
\]
\end{lemma}
\begin{proof}
We compute
\[
f_\sigma(a) = \int_a^b \phi_\sigma(a-y)\,dy = \int_a^{a+\sigma} \phi_\sigma(a-y)\,dy = \int_{-\sigma}^0\phi_\sigma(y)\,dy = \frac12
\]
and by symmetry $f_\sigma(b)=\frac12$.

Furthermore for $x\in (a-\sigma, b+\sigma)$ we compute
\[
f_\sigma'(x) = \int_a^b \phi_\sigma'(x-y)\,dy = \int_{x-b}^0 \phi_\sigma'(y)\,dy + \int_0^{x-a} \phi_\sigma'(y)\,dy.
\]
The first term on the right is nonnegative if $x \leq b$ and nonpositive if $x\geq b$ and the second term is nonpositive if $x\geq a$ and nonnegative if $x\leq a$. By symmetry of $\phi_\sigma$ we then conclude that $f_\sigma'(x) \geq 0$ if $x\leq \frac{a+b}2$ and $f_\sigma'(x) \leq 0$ if $x\geq \frac{a+b}2$. Moreover we have $f_\sigma'(a) >0$ and $f_\sigma'(b)<0$. 
\end{proof}


\begin{lemma}\label{lem:lessthanahalf}
Let $\phi_\sigma \in \mathcal{K}$, $z\in \B_\omega$, and $\sigma \leq \omega$, then for every $x\in \R\setminus \supp z$ we have $\phi_\sigma*z(x) < \frac12$.
\end{lemma}
\begin{proof}
Let $x \in \R\setminus \supp z$, then there exist $a<b$ such that $b-a \geq \omega$, $x \in (a,b)$, and $z(y)=0$ for all $y\in (a,b)$. Define $z_0 := \chi_{(-\infty, a)} + \chi_{(b, \infty)}$, then
\[
\phi_\sigma * z(x) \leq \phi_\sigma * z_0(x) 
= \int_\R \phi_\sigma(x-y)\, dy - \int_a^b \phi_\sigma(x-y)\,dy
= 1- \phi_\sigma*\chi_{[a,b]}(x)
< \frac12.
\]
The final inequality follows since $\sigma \leq \omega \leq b-a$ and thus by Lemma~\ref{lem:onebarconvolution} $\phi_\sigma*\chi_{[a,b]} > \frac12$ on $(a, b)$.
\end{proof}

The following lemma allows us to consider only minimisers of $F_1$ that have less interfaces than $z$ or are equal to $z$.
\begin{lemma}\label{lem:firstvar-F1}
Let $\phi_\sigma \in \mathcal{K}$, $z\in \B_{\omega}$ 
and let $u$ be a minimiser of $F_1$ over $\B$. Denote by $x_i$ the locations of the interfaces of $u$, then we have for every $i$
\[
f_\sigma(x_i)=\frac12.
\]
Consequently if $\sigma \leq \omega$, $x_i$ is the location of an interface of $z$ for every $i$.
\end{lemma}
\begin{proof}
Assume without loss of generality that $z\neq 0$. Let $u$ minimise  $F_1$ over $\B$. 
We show that vanishing first variation of $u$ implies that at any interface $x_i$, we must have $f_\sigma (x_i) = \frac12$.  To this end, consider an interface of transition from $u=1$ to $u=0$ at $x_i$ (the other case is treated similarly). By considering a perturbation consisting of extending the $u = 1$ bar up to   $x_i + t$ for $t$ small, one obtains 
no change in the total variation and a change in the fidelity term of 
\[   \int_{x_i}^{x_i + t} \Bigl((1 - f_\sigma)^2  -  f_\sigma^2\Bigr)\, dx \, = \, 
   \int_{x_i}^{x_i + t}  (1 - 2f_\sigma) \, dx. \]
Differentiating with respect to $t$ and setting $t = 0$ gives $1 - 2 f_\sigma (x_i) = 0$. 

Let now $\sigma \leq \omega$. By Lemma~\ref{lem:onebarconvolution} if $z$ consists of one bar only the $\frac12$-level set of $f_\sigma$ is exactly the set of locations of the interfaces of $z$. If $z$ has more bars Lemma~\ref{lem:lessthanahalf} assures that the $\frac12$-lower level set is not affected by the convolutions of different bars interacting.



\end{proof}

\bigskip

Lemma~\ref{lem:firstvar-F1} tells us that, if $\sigma\leq \omega$, any candidate for minimising $F_1$ over $\B$ not equal to $z$, should have less interfaces than $z$ which are located at places where $z$ also has an interface. We use this to complete the proof of Theorem~\ref{thm:minimofF2}. First, we introduce a notation that will be used frequently in what follows.
For $\sigma>0$, $a<b$, and $x\in\R$ define the functions
\begin{equation}\label{eq:Ipm}
I_\pm^\sigma(x, a, b) := \frac1\sigma \int_a^b \left(1\pm \frac{x-y}\sigma\right)\, dy = \frac1\sigma \left( (b-a) \left( 1\pm\frac{x}\sigma\right)  \mp \frac{b^2-a^2}{2\sigma} \right).
\end{equation}

Note that this definition is tailored to the needs of the hat function $\hat \phi_\sigma$. If we want to reproduce the calculations that follow for a general convolution kernel $\phi_\sigma\in \mathcal{K}$ we can write $\phi_\sigma$ as in (\ref{eq:Ksigma}) and define $I_\pm^\sigma$ as
\[
I_+^\sigma(x, a, b) := \int_a^b p(y-x, \sigma)\, dy \quad \text{and} \quad I_-^\sigma(x, a, x) := \int_a^b p(x-y, \sigma)\, dy.
\]

We will only state the results for the hat function and hence use the definitions in (\ref{eq:Ipm}).

\bigskip

\noindent{\bf Proof of Theorem~\ref{thm:minimofF2}, part~\ref{item:F1zminoverB}.}
Let $\B \ni u_0\neq z$ be a minimiser of $F_1$. 
By Lemma \ref{lem:firstvar-F1}, the number of interfaces of $u_0$ must be less than the number of interfaces of $z$ and the location of every interface of $u_0$ coincides with the location of an interface of $z$. 
Therefore there exists a connected interval $N\subset[0,1]$ such that  $|N| \ge \omega$ and either
\begin{itemize}
\item $z=0$ and $u_0=1$ on $N$, or
\item $z=1$ and $u_0=0$ on $N$.
\end{itemize}
First assume the former case,  and let
\[
\hat u := \left\{ \begin{array}{ll} u_0 & \text{on } N^c \\ u_0-1=z & \text{on } N. \end{array}\right.
\]
We compute
\begin{equation}\label{eq:jumpestimate}
\int_{\R} |\hat u'| \leq 2+\int_{\R} |u_0'|
\end{equation}
and
\begin{align}
\|u_0 - f_\sigma\|_{L^2(\R)}^2 &= \|u_0 - \hat u + \hat u - f_\sigma\|_{L^2(\R)}^2\notag\\
&= \|u_0-\hat u\|_{L^2(\R)}^2 + \|\hat u - f_\sigma\|_{L^2(\R)}^2 + 2 \int_{\R} (u_0-\hat u) (\hat u - f_\sigma)\notag\\
&= \|\hat u - f_\sigma\|_{L^2(\R)}^2 + |N| + 2 \int_{\R} (u_0-\hat u) (\hat u - f_{\sigma}),\label{eq:L2normestimate}
\end{align}
where we have used that
\[
\|u_0-\hat u\|_{L^2(\R)}^2 = |N|.
\]


Let $a\in [0,1]$ be such that $N=[a, a+|N|]$ and
thus $N^c\cap[0,1]=[0,a)\cup(a+|N|,1]$, then we compute
\begin{align}\label{eq:twointegrals}
-2 \int_{\R} (u_0-\hat u) (\hat u - f_\sigma) &= 2 \int_N \int_{\R} \hat\phi_{\sigma}(x-y) z(y) \, dy \, dx 
\leq	 2 \int_N \int_{N^c\cap[0,1]} \hat\phi_{\sigma}(x-y)\, dy \, dx \nonumber\\
&= 2 \int_a^{a+|N|} \left\{\int_0^a \hat\phi_{\sigma}(x-y)\,dy + \int_{a+|N|}^1 \hat\phi_{\sigma}(x-y)\, dy\right\}\, dx.
\end{align}

Integrals as those in the right hand side of (\ref{eq:twointegrals}) are commonplace in the proofs of this paper. It is therefore very illustrative to work out one of them in detail. Let us consider
\[
\int_a^{a+|N|} \int_0^a \hat\phi_{\sigma}(x-y)\,dy \,dx.
\]
Per definition $\hat\phi_\sigma(x-y)$ is zero if $|x-y|\geq\sigma$ and on its support its value is given by $\frac1\sigma \left(1+\frac{x-y}\sigma\right)$ if $x-y\in (-\sigma, 0)$ and by $\frac1\sigma \left(1-\frac{x-y}\sigma\right)$ if $x-y\in (0,\sigma)$. Let us fix $x\in[a, a+|N|]$ for the moment and remember that $y\in(0,a)$ in the integral, then $\phi_\sigma(x-y) = \frac1\sigma \left(1-\frac{x-y}\sigma\right)$ if
\[
y \in (x-\sigma, x) \cap (0, a) = \left\{ \begin{array}{ll} 
\emptyset & \text{if } x-\sigma < x < 0 < a,\\
(0, x) & \text{if } x-\sigma < 0 < x < a,\\
(0, a) & \text{if } x-\sigma < 0 < a < x,\\
(x-\sigma, x) & \text{if } 0 < x-\sigma < x < a,\\
(x-\sigma, a) & \text{if } 0< x-\sigma< a < x,\\
\emptyset & \text{if } 0 < a < x-\sigma < x.
\end{array}\right.
\]
Because $x\in[a, a+|N|]$ we can rule out some of these cases\footnote{For many of the similar calculations in the rest of this paper, $x\in \R$ and this kind of simplification will not be possible.} and end up with
\[
y \in \left\{ \begin{array}{ll}
(0, a) & \text{if } x \in (-\infty, \sigma) \cap (a, \infty) = \left\{ \begin{array}{ll} \emptyset & \text{if } \sigma \leq a,\\ (a, \sigma) & \text{if } a < \sigma,\end{array}\right.\\
(x-\sigma, a) &\text{if } x \in (\sigma, a+\sigma) \cap (a, \infty) = \left\{ \begin{array}{ll} (a, a+\sigma) & \text{if } \sigma \leq a,\\ (\sigma, a+\sigma) & \text{if } a < \sigma,\end{array}\right.\\
\emptyset & \text{if } x\in (a+\sigma, a+|N|).
\end{array}\right.
\]
We see that we have to distinguish between the cases $\sigma \leq a$ and $a < \sigma$. Similarly we find that $\hat\phi_\sigma(x-y) = \frac1\sigma \left(1+\frac{x-y}\sigma\right)$ if $y \in (x, x+\sigma) \cap (0, a)$, which is the empty set because of the restrictions on $x$.

This now leads us to the computation
\[
\int_a^{a+|N|} \int_0^a \hat\phi_{\sigma}(x-y)\,dy \,dx = \left\{ \begin{array}{ll}
\frac1\sigma \int_a^{a+\sigma} \int_{x-\sigma}^a \left(1-\frac{x-y}\sigma\right)\,dy\,dx & \text{if } \sigma \leq a,\\
\frac1\sigma \left[ \int_a^\sigma \int_0^a \left(1-\frac{x-y}\sigma\right)\,dy\,dx + \int_\sigma^{a+\sigma} \int_{x-\sigma}^a \left(1-\frac{x-y}\sigma\right)\,dy\,dx\right] & \text{if } a < \sigma.
\end{array}\right.
\]
Because all the integrands are positive, in the case $a<\sigma$ we can estimate
\[
\int_a^\sigma \int_0^a \left(1-\frac{x-y}\sigma\right)\,dy\,dx \leq \int_a^\sigma \int_{x-\sigma}^a \left(1-\frac{x-y}\sigma\right)\,dy\,dx.
\]
Therefore we conclude that for both $\sigma \leq a$ and $a<\sigma$
\[
\int_a^{a+|N|} \int_0^a \hat\phi_{\sigma}(x-y)\,dy \,dx \leq \int_a^{a+\sigma} I_-^\sigma(x, x-\sigma, a)\, dx.
\]

In a similar fashion we compute
\[
\int_a^{a+|N|}  \int_{a+|N|}^1 \hat\phi_{\sigma}(x-y)\, dy\, dx \leq \int_{a+|N|-\sigma}^{a+|N|} I_+^\sigma(x, a+|N|, x+\sigma)\,dx.
\]
During this computation we need to distinguish between the cases $a+|N| \leq 1-\sigma$ and $a+|N| > 1-\sigma$, but as before this distinction doesn't play a role in the final estimate.

Continuing from (\ref{eq:twointegrals}) we now find
\begin{align*}
-2 \int_{\R} (u_0-\hat u) (\hat u - \hat\phi_{\sigma} * z) &\leq 2 \left\{ \int_a^{a+\sigma} I_-^\sigma(x, x-\sigma, a)\, dx + \int_{a+|N|-\sigma}^{a+|N|} I_+^\sigma(x, a+|N|, x+\sigma)\,dx\right\}\\ 
&= \frac23 \sigma.
\end{align*}
Using this in (\ref{eq:jumpestimate}--\ref{eq:L2normestimate}) together with $|N| \geq \omega$ we find
\[
F_1 (\hat u) \leq F_1 (u_0) + 2 + \lambda \left(\frac23 \sigma - \omega \right) < F_1(u_0),
\]
where the second inequality follows by condition (\ref{eq:conditionsigmalambda}). This contradicts $u_0$ being a minimiser of $F_1$.

Next we consider the second case, i.e. $z=1$ and $u_0=0$ on $N$. We define
\[
\bar u := \left\{ \begin{array}{ll} u_0 & \text{on } N^c \\ u_0+1=z & \text{on } N. \end{array}\right.
\]
As in the first case,
we will find an estimate for the integral in the brackets in the right hand side of (\ref{eq:L2normestimate}), but for $\bar u$ instead of $\hat u$:
\[
-2 \int_{\R} (u_0-\bar u) (\bar u - \hat\phi_{\sigma} * z) = 2 \int_N (1-\hat\phi_{\sigma} * z) = 2 |N| - 2 \int_N \hat\phi_{\sigma} * z.
\]

Again we write $N=[a, a+|N|]$ and we compute
\begin{align*}
\int_N \phi_{\sigma} * z &= \int_N \int_0^1 \hat\phi_{\sigma}(x-y) z(y) \, dy \,dx 
\geq \int_N \int_N \hat\phi_{\sigma}(x-y) z(y)\,dy\,dx\\
&= \int_a^{a+|N|}\int_a^{a+|N|} \hat\phi_{\sigma}(x-y)\,dy\,dx\\
&= \int_a^{a+\sigma} I_-^\sigma(x, a, x)\, dx + \int_{a+\sigma}^{a+|N|} I_-^\sigma(x, x-\sigma, x) \,dx\\
&\hspace{0.4cm}  + \int_a^{a+|N|-\sigma} I_+^\sigma(x, x, x+\sigma) \,dx + \int_{a+|N|-\sigma}^{a+|N|} I_+^\sigma(x, x, a+|N|)\,dx\\
&= |N| - \frac13\sigma.
\end{align*}
As in the first case we now find
\[
F_1(\bar u) \leq F_1(u_0) + 2 + \lambda \left(\frac23 \sigma - \omega \right) < F_1 (u_0),
\]
which is again a contradiction with $u_0$ being a minimiser.
Therefore the only candidate for a minimiser is $u=z$ and hence by Lemma~\ref{lemma0} $u=z$ is the unique minimiser. 
\qed

\begin{remark}
In the above proof everything up to and including (\ref{eq:twointegrals}) is independent of the choice of specific blurring kernel and we could have used any $\phi_\sigma \in \mathcal{K}$. The explicit calculations that follow in the remainder of the proof depend on our choice $\phi_\sigma=\hat\phi_\sigma$, but can be redone for a different choice of kernel as explained in the paragraphs preceding the proof.
\end{remark}

\subsection{Proofs of Theorem~\ref{thm:F4}, part~\ref{item:F3zminoverB} and Corollary~\ref{thm:F3}, part~\ref{item:F2zminoverB}}

We now turn to $F_3$. 

\begin{lemma}\label{lemma1} 
Let  $z_1, z_2 \in \B$. Both $z_1$ and $z_2$ on $[0,1]$ consist of a finite collection of subintervals of $[0,1]$, i.e. alternating bars and spaces. 
Let $t_i$, $i = 1 \dots n$ and  $t_i'$, $i = 1 \dots n'$ denote the right hand sides of the intervals of 
 $z_1$ and $z_2$ respectively.  In particular $t_n=t_{n'}'=1$.
If $n > n'$,  then there  exists an 
interval  $N\subset [0,1]$ such that $[t_i, t_{i+1}] \subset N$ 
for some $i$; and for all $ x \in N$, either
\begin{equation}\label{def-N}
 z_1(x) =0  \,\,\, {\rm and} \,\,\, z_2 (x) =1  \qquad {\rm  or} \qquad z_1 (x) =1 \,\,\, {\rm  and} \,\,\, z_2 (x) =0. 
 \end{equation}
In particular, if $z_1\in \B_{\omega}$, then $|N|\geq \omega$.
\end{lemma}

\begin{proof}
First assume that $z_1$ starts with a bar and $z_2$ starts with a space, i.e. $z_1=1$ on $[0, t_1]$ and $z_2=0$ on $[0, t_1']$. If $t_1\leq t_1'$ then $[0, t_1]\subset N$. Suppose $t_1' \leq t_1$. If the conclusion of the lemma is false, then for all $i \leq n'$, $t_i' < t_i$. This is a contradiction since $t'_{n'} = t_n = 1$. 
If $z_1$ starts with a space and $z_2$ starts with a bar we arrive at a similar conclusion.

Now assume that $z_1$ and $z_2$ both start with a bar (the situation in which both start with a space is similar). Note that $z_1 = 1$ on $[0, t_1]$ and $z_2 = 1$ on $[0, t_1']$. 
Suppose $t_1 \le t_1'$. Then if the conclusion of the lemma is false, we must have $t_i' < t_{i+1}$ for  $i=1 \dots n'$ which implies  $1 = t'_{n'} < 1$.   
Suppose  $t_1 > t_1'$. If for some $i>1$, we have $t_i' \ge t_i$, then the previous argument again gives a contradiction. Thus we must have $t_i > t_i'$ for all $i = 2 \dots n' -1$. But then 
(\ref{def-N}) must hold on one of the intervals $[t_i, t_{i+1}]$, for $i \ge n'$. 
\end{proof}

\begin{lemma}\label{lem:F3partresult}
Let $z
\in \B_{\omega}^{ij}$ for some $i, j\in \{0,1\}$, $\rho, \sigma \leq \frac{\omega}2$, $\phi_\sigma = \hat\phi_\sigma$ and define 
\begin{equation}\label{eq:thisfunctionofrhoandsigma}
f(\rho, \sigma) := \left\{ \begin{array}{ll} \frac1{\rho^2} \Bigl(-\sigma^3 + 5 \rho \sigma^2 + 10 \rho^3\Bigr) & \text{if } \sigma\leq \rho,\\ \frac1{\sigma^2} \Bigl( -\rho^3+5 \sigma \rho^2 + 10 \sigma^3\Bigr) & \text{if } \rho\leq\sigma.\end{array} \right.
\end{equation}
Let $\lambda$, $\rho$, and $\sigma$ satisfy in addition
\begin{equation}\label{eq:conditiononrhosigmalambda2}
\frac2\lambda + \frac1{15} \Bigl(7 \rho + f(\rho, \sigma)\Bigr) < \omega.
\end{equation}
If $u \in \B^{ij}$ is a minimiser of $F_3$ over $\B^{ij}$, then
\[ \int_\R |u'| \, \geq \, \int_\R |z'|.\]
\end{lemma}

\bigskip


\begin{proof}[Proof of Lemma~\ref{lem:F3partresult}]

We prove this by contradiction. Let $\B^{ij}\ni u_0 \neq z$ be a minimiser of $F_3$ in $\B^{ij}$ and assume that $u_0$ has less interfaces than $z$, i.e. $\int_\R |u_0'| \, < \, \int_\R |z'|$.
By Lemma \ref{lemma1},  there exists a connected interval $N\subset[0,1]$ such that  $|N| \ge \omega$ and either
\begin{itemize}
\item $z=0$ and $u_0=1$ on $N$, or
\item $z=1$ and $u_0=0$ on $N$.
\end{itemize}
Define
\[
\hat u := \left\{\begin{array}{ll} u_0 & \text{on } N^c,\\ z & \text{on } N,  \end{array}\right.
\]
then
\[
\int |\hat u'| \leq \int |u_0'| + 2
\]
and
\begin{align*}
\|\hat\phi_\rho * u_0 - \hat\phi_\sigma*z\|_{L^2(\R)}^2 &= \|\hat\phi_\rho * \hat u -\hat\phi_\sigma*z\|_{L^2(\R)}^2 + \|\hat\phi_\rho * (u_0-\hat u)\|_{L^2(\R)}^2\\ &\hspace{0.4cm}+ 2 \int_\R \Big(\hat\phi_\rho * (u_0-\hat u)\Bigr) \cdot \Big(\hat\phi_\rho * \hat u - \hat\phi_\sigma* z\Big),
\end{align*}
from which we conclude that
\begin{equation}\label{eq:F4inequality}
F_3(\hat u) \leq  F_3(u_0) + 2 - \lambda \Biggl( \|\hat\phi_{\rho} * (u_0-\hat u)\|_{L^2}^2 + 2 \int_\R \Big(\hat\phi_{\rho} * (u_0-\hat u)\Big) \cdot \Bigl(\hat\phi_{\rho}*\hat u - \hat\phi_{\sigma}*z\Bigr)\Biggr).
\end{equation}
Because
$u_0-\hat u = \pm \chi_N$, 
Lemma~\ref{lem:fsigmasquared} gives
\[
 \|\hat\phi_{\rho} * (u_0-\hat u)\|_{L^2}^2 = |N| - \frac7{15} \rho.
\]
Next we again distinguish two cases: Case I in which $u_0=1$ and $\hat u = z = 0$ on $N$ and Case II in which $u_0=0$ and $\hat u = z = 1$ on $N$. We first treat Case I:
\begin{align*}
&\hspace{0.4cm}2 \int_\R \hat\phi_\rho * (u_0-\hat u) \Bigl(\hat\phi_\rho * \hat u - \hat\phi_\sigma * z \Bigr)\\
&= 2 \int_R \int_N \hat\phi_\rho(x-y)\,dy \biggl(\int_{N^c\cap[0,1]} \hat\phi_\rho(x-w) \hat u(w)\,dw - \int_{N^c\cap[0,1]} \hat\phi_\sigma(x-w) z(w)\,dw\biggr)\,dx\\
&\geq -2 \int_\R \int_N \hat\phi_\rho(x-y)\,dy \int_{N^c\cap[0,1]} \hat\phi_\sigma(x-w)\,dw\,dx.
\end{align*}
Now we subdivide Case I into two subclasses: Case Ia in which $\sigma\leq\rho$ and Case Ib in which $\rho\leq\sigma$. For Case Ia we compute
\[
-2 \int_\R \int_N \hat\phi_\rho(x-y)\,dy \int_{N^c\cap[0,1]} \hat\phi_\sigma(x-w)\,dw\,dx = \frac1{15 \rho^2} \Bigl(\sigma^3 - 5 \rho \sigma^2 - 10 \rho^3\Bigr).
\]
For details of this computation we refer to (\ref{eq:lemCaseIa}) in Appendix~\ref{sec:moredetails}.


In Case Ib the computation is
\[
-2 \int_\R \int_N \hat\phi_\rho(x-y)\,dy \int_{N^c\cap[0,1]} \hat\phi_\sigma(x-w)\,dw\,dx = \frac1{15\sigma^2} \Bigl(\rho^3 - 5 \sigma \rho^2 - 10 \sigma^3\Bigr),
\]
the details of which can be found in (\ref{eq:lemCaseIb}) in Appendix~\ref{sec:moredetails}.

In Case II we compute
\begin{align*}
&\hspace{0.4cm}2 \int_\R \hat\phi_\rho * (u_0-\hat u) \Bigl(\hat\phi_\rho * \hat u - \hat\phi_\sigma * z \Bigr)\\
&= -2\int_\R \int_N \hat\phi_\rho(x-y)\,dy \biggl(\int_\R \hat\phi_\rho(x-w)\hat u(w)\,dw - \int_\R \hat\phi_\sigma(x-w) z(w)\,dw\biggr)\,dx\\
&\geq -2\int_\R \int_N \hat\phi_\rho(x-y)\,dy \biggl(\int_\R \hat\phi_\rho(x-w)\,dw - \int_N\hat\phi_\sigma(x-w)\,dw\biggr)\,dx.
\end{align*}

For the first term we find
\[
-2\int_\R \int_N \hat\phi_\rho(x-y)\,dy \int_\R \hat\phi_\rho(u-w)\,dw\,dx =-2N.
\]
Details of this calculation are given in (\ref{eq:lemCaseIIfirstterm}) in Appendix~\ref{sec:moredetails}.

For the second term again we need to subdivide into Case IIa in which $\sigma\leq\rho$ and Case IIb in which $\rho\leq\sigma$. For Case IIa we compute
\[
2 \int_\R \int_N \hat\phi_\rho(x-y) \,dy \int_N \hat\phi_\sigma(x-w)\,dw\,dx = 2N + \frac1{15 \rho^2} \Bigl(\sigma^3 - 5 \rho \sigma^2 - 10 \rho^3\Bigr).
\]
For more details of this computation see (\ref{eq:lemCaseIIa}) in Appendix~\ref{sec:moredetails}.

In Case IIb we can repeat the calculation with $\rho$ and $\sigma$ interchanged to get
\[
2 \int_\R \int_N \hat\phi_\rho(x-y) \,dy \int_N \hat\phi_\sigma(x-w)\,dw\,dx = 2N + \frac1{15 \sigma^2} \Bigl(\rho^3 - 5 \rho \sigma^2 - 10 \sigma^3\Bigr).
\]

Using the combined results of Cases I and II in inequality (\ref{eq:F4inequality}) leads to
\begin{align*}
F_3(\hat u) &\leq  F_3(u_0) + 2 - \lambda \biggl( |N| -\frac1{15} \Big(7 \rho + f(\rho, \sigma)\Big) \biggr)\\
&\leq  F_3(u_0) + 2 - \lambda \biggl( \omega -\frac1{15} \Big(7 \rho + f(\rho, \sigma)\Big) \biggr)\\
&< F_3(u_0),
\end{align*}
where the final inequality follows from (\ref{eq:thisfunctionofrhoandsigma}) - (\ref{eq:conditiononrhosigmalambda2}). This contradicts the fact that $u_0$ is a minimiser of $F_3$. 

\end{proof}

\begin{remark}\label{rem:roleofconditions}
In the proof of Lemma~\ref{lem:F3partresult} we have used the conditions $z\in \B_\omega^{ij}$, $u_0\in \B^{ij}$, and $\rho, \sigma \leq \frac\omega2$ but it might not be immediately clear where. They allow us to order the endpoints of the intervals of integration that occur in the integrals in Appendix~\ref{sec:moredetails}. In particular $z\in \B_\omega^{ij}$ and $u_0\in \B^{ij}$ imply that the interval $N$ on which $z$ and $u_0$ differ is located at least a distance $\omega$ away from the endpoints of the interval $[0, 1]$, i.e. $a \geq \omega$ and $a+|N| \leq 1-\omega$. If we also take into account the conditions $\rho, \sigma \leq \frac\omega2$ we have the ordering, for $\sigma\leq\rho$,
\begin{align*}
&-\rho \leq -\sigma \leq 0 \leq a-\rho-\sigma \leq a-\sigma \leq a \leq a+\sigma \leq a+\rho \leq a+|N|-\rho \leq a+|N|-\sigma\\
&\hspace{0.65cm}\leq a+|N| \leq a+|N|+\sigma \leq a+|N|+\rho \leq 1-\rho \leq 1-\sigma \leq 1 \leq 1+\sigma \leq 1+\rho
\end{align*}
and an analogous one for $\rho\leq \sigma$. These orderings are important when determining exactly which $I_\pm^\sigma(x, a, b)$ contribute over which $x$-intervals to integrals like
\[
\int_\R \int_N \hat\phi_\rho(x-y)\,dy \int_{N^c\cap[0,1]} \hat\phi_\sigma(x-w)\,dw\,dx.
\]

Loosening the condition $z \in \B_\omega^{ij}$ to $z \in \B_\omega$ and consequently $u_0\in \B^{ij}$ to $u_0\in \B$ is possible in principle, but will give rise to more possible orderings of the kind above and separate calculations of all the integrals involved need to be done for each possible ordering. It is not expected however that this will influence the end result by much if at all. 

\end{remark}

\begin{remark}
Up to and including (\ref{eq:F4inequality}) the steps in the proof of Lemma~\ref{lem:F3partresult} are independent of the specific choice of kernels $\phi_\sigma$ and $\phi_\rho$, but the calculations that make up the remainder of the proof do depend on the explicit choice $\phi_\sigma=\hat\phi_\sigma$. In order to derive similar results for other kernels we need to redo those computations with an explicitly given alternative choice.
\end{remark}




\bigskip

The result for $F_2$ in Corollary~\ref{thm:F3}, part~\ref{item:F2zminoverB} follows as a direct consequence of Theorem~\ref{thm:F4}, part~\ref{item:F3zminoverB} for $F_3$ by choosing $\rho=\sigma$. However, the fact that the fidelity term in $F_2$ vanishes if and only if $u=z$ allows for a direct proof as well.
\\

\noindent{\bf Proof of Corollary~\ref{thm:F3}, part~\ref{item:F2zminoverB}}:
Since for $F_2$, the fidelity term vanishes at $u = z$  any  potential competitor must have strictly less interfaces than $z$. The result follows then immediately from Lemma~\ref{lem:F3partresult} with $\rho=\sigma$ and Lemma~\ref{lemma0}.
\qed

\bigskip

\bigskip

To complete the proof of Theorem~\ref{thm:F4} we need a result that tells us that, if $\sigma \leq \rho$, a minimiser of $F_3$ is either equal to $z$ or has strictly less interfaces. Lemma~\ref{lem:F3firstvar} will provide exactly this. First we need some preparatory lemmas.

\begin{lemma}\label{lem:levelhalfrhosigmaz}
Let $z\in \B_\omega$, $\rho, \sigma \leq \frac\omega2$, and $\phi_\sigma \in \mathcal{K}\cap C(\R)$, then the level-$\frac12$ set of $\phi_\rho*f_\sigma$ consists of exactly the locations of the interfaces of $z$. Furthermore the upper level-$\frac12$ set where $\phi_\rho*\phi_\sigma*z \geq \frac12$ is $\supp z$.
\end{lemma}
\begin{proof}
If $z=0$ the results follow trivially. We assume now $z\neq 0$. First we consider the case of a bar code with only one bar. Let $a<b$ be such that $b-a \geq \omega$ and define $z:= \chi_{[a, b]}$.
Since the convolution of two symmetric unimodal functions is again a symmetric unimodal function (see \cite{Uhrin84, EatonPerlman91} and references therein) we find
%
that $\phi_\sigma*z$ is a unimodal function with mode at $x_0:=\frac{a+b}2$, i.e. $\phi_\sigma*z$ is non-decreasing for $x\geq x_0$ and non-increasing for $x\leq x_0$, and symmetric around $x=x_0$. Since $\phi_\rho$ is unimodal with mode at $x=0$ and symmetric around $x=0$ we 
conclude that $\phi_\rho*\phi_\sigma*z$ is unimodal with mode at $x=x_0$ and symmetric around $x=x_0$. Therefore for all $x \leq a$ and all $x\geq b$
\begin{equation}\label{eq:nonstrictineq1}
\phi_\rho*\phi_\sigma*z(x) \leq \phi_\rho*\phi_\sigma*z(a)=\phi_\rho*\phi_\sigma*z(b)
\end{equation}
and for all $x\in [a,b]$
\begin{equation}\label{eq:nonstrictineq2}
\phi_\rho*\phi_\sigma*z(x) \geq \phi_\rho*\phi_\sigma*z(a)=\phi_\rho*\phi_\sigma*z(b).
\end{equation}
In the sense of distributions we have
\[
z' = \delta_a-\delta_b
\]
where $\delta_x$ is the Dirac delta measure at $x$. Hence
\[
\phi_\rho*\phi_\sigma*z'(x) = \phi_\rho*\phi_\sigma(x-a) - \phi_\rho*\phi_\sigma(x-b).
\]
Because $\phi_\rho*\phi_\sigma$ is unimodal with maximum at $0$ we deduce
that $(\phi_\rho * \phi_\sigma * z)'(a)>0$ and $(\phi_\rho * \phi_\sigma * z)'(b)<0$. Combined with (\ref{eq:nonstrictineq1}) and (\ref{eq:nonstrictineq2}) this implies that for all $x < a$ and all $x > b$
\[
\phi_\rho*\phi_\sigma*z(x) < \phi_\rho*\phi_\sigma*z(a)=\phi_\rho*\phi_\sigma*z(b)
\]
and for all $x\in (a,b)$
\[
\phi_\rho*\phi_\sigma*z(x) > \phi_\rho*\phi_\sigma*z(a)=\phi_\rho*\phi_\sigma*z(b).
\]
We now explicitly compute the value $\phi_\rho*\phi_\sigma*z(a)$. 
\begin{align}
\phi_\rho*\phi_\sigma*z(a) &= \int_\R \int_\R \phi_\rho(a-x) \phi_\sigma(x-y) \chi_{[a,b]}(y)\, dy\, dx\notag\\
&= \int_{a-\rho}^{a+\rho} \int_a^{x+\sigma} \phi_\rho(a-x) \phi_\sigma(x-y)\, dy\, dx\notag 
= \int_{-\rho}^\rho \int_{-\sigma}^{-z} \phi_\rho(z) \phi_\sigma(q)\, dq\, dz\notag\\
&= \int_{-\rho}^\rho \int_{-\sigma}^{0} \phi_\rho(z) \phi_\sigma(q)\, dq\, dz - \int_{-\rho}^\rho \int_{-z}^0 \phi_\rho(z) \phi_\sigma(q)\, dq\, dz\notag\\
&= \frac12 -  \int_{-\rho}^\rho \int_{-z}^0 \phi_\rho(z) \phi_\sigma(q)\, dq\, dz.\label{eq:1/2}
\end{align}
In the third equality we have used the change of variables
\[
\left(\begin{array}{c} z\\ q\end{array}\right) = \left(\begin{array}{c} a\\ 0\end{array}\right) + \left(\begin{array}{cc} -1&0\\ 1&-1\end{array}\right)\left(\begin{array}{c} x\\ y\end{array}\right).
\]
The last equality follows by symmetry of $\phi_\sigma$ and the fact that $\phi_\rho$ and $\phi_\sigma$ have unit mass.

Because
\begin{align*}
\int_{-\rho}^0 \int_{-z}^0 \phi_\rho(z) \phi_\sigma(q)\, dq \, dz &= \int_{\rho}^0 \int_z^0 \phi_\rho(-z) \phi_\sigma(q)\, dq \, d(-z) = -\int_0^\rho \int_0^z \phi_\rho(z) \phi_\sigma(q)\, dq \, dz\\
&= -\int_0^\rho \int_{-z}^0 \phi_\rho(z) \phi_\sigma(q)\, dq \, dz
\end{align*}
we have
\[
\int_{-\rho}^\rho \int_{-z}^0 \phi_\rho(z) \phi_\sigma(q)\, dq\, dz = \int_{-\rho}^0 \int_{-z}^0 \phi_\rho(z) \phi_\sigma(q)\, dq\, dz + \int_0^\rho \int_{-z}^0 \phi_\rho(z) \phi_\sigma(q)\, dq\, dz
= 0
\]
and hence by (\ref{eq:1/2})
\begin{equation}\label{eq:equaltoahalf}
\phi_\rho*\phi_\sigma*z(a) = \frac12.
\end{equation}

This proves the result if $z$ has only one bar. If $z$ has more bars then we prove that the $\frac12$-lower level set of $\phi_\rho*f_\sigma$ is the same as the $\frac12$-lower level set of $f_\sigma$ in a similar fashion as the $\frac12$-lower level set was identified in the proof of Lemma~\ref{lem:lessthanahalf}. Let $x\in \R\setminus \supp z$, then there exist $c<d$ such that $d-c \geq \omega$, $x\in (c, d)$, and $z(y)=0$ for all $y\in (c, d)$. Define $z_0 := \chi_{(-\infty, c)} + \chi_{(d, \infty)}$, then
\begin{align*}
\phi_\sigma * f_\sigma(x) &\leq \phi_\rho*\phi_\sigma * z_0(x) 
= \int_\R \phi_\rho*\phi_\sigma(x-y)\, dy - \int_c^d \phi_\rho*\phi_\sigma(x-y)\,dy\\
&= 1- \phi_\rho*\phi_\sigma*\chi_{[c,d]}(x)
< \frac12.
\end{align*}
The last inequality follows from $\phi_\rho*\phi_\sigma*\chi_{[c, d]} > \frac12$ on $(c, d)$ as proven above.

\end{proof}


\begin{lemma}\label{lem:decreasing}
Let $z\in \B_\omega$, $\phi_\rho\in \mathcal{K}_3$ and $\rho \leq \frac\omega2$. Fix $x\in \supp z$, then the function
\[
(0, \rho] \to \R: \tau \mapsto \phi_\tau * \phi_\rho*z(x)
\]
is non-increasing.
\end{lemma}
\begin{proof}
First assume that $z =\chi_{[a, b]}$ for some $a<b$ satisfying $b-a \geq \omega$. We compute
\[
\frac{\partial}{\partial \tau} \phi_\tau*\phi_\rho*z(x) = \int_\R \int_a^b \phi_\rho(y-w) \frac{\partial}{\partial \tau} \phi_\tau(x-y) \, dw\, dy.
\]
As in (\ref{eq:Ksigma}) we write $
\phi_\tau(x) =  p(-x, \tau) \chi_{[-\tau, 0]}(x) + p(x, \tau) \chi_{[0, \tau]}(x)$ and thus, by continuity of $\phi_\tau$ in $\tau$,

\[
\frac{\partial}{\partial \tau} \phi_\tau(x) = \chi_{[-\tau, 0]}(x) \frac{\partial}{\partial \tau} p(-x, \tau) +  \chi_{[0, \tau]}(x)\frac{\partial}{\partial \tau}p(x, \tau)
\]
and thus
\begin{align*}
\frac{\partial}{\partial \tau} \phi_\tau*\phi_\rho*z(x) &= \int_x^{x+\tau} \int_a^b \phi_\rho(y-w) \frac{\partial}{\partial \tau} p(y-x, \tau)\, dw\, dy\\ &\hspace{0.7cm}+  \int_{x-\tau}^x \int_a^b \phi_\rho(y-w) \frac{\partial}{\partial \tau} p(x-y, \tau)\, dw\, dy.
\end{align*}
Using the substitution of variables $\displaystyle \left(\begin{array}{c} \tilde y\\ \tilde w\end{array}\right) = \left(\begin{array}{c} -x\\ x\end{array}\right) + \left(\begin{array}{cc} 1&0\\ 0&-1\end{array}\right) \left(\begin{array}{c} y\\ w\end{array}\right)$, using the symmetry of $\phi_\sigma$, then writing $\tilde x = x-a$ and $c=b-a$ and finally dropping the tildes, allows us to rewrite the integrals above as the integral in (\ref{eq:condpsit}) with $\rho$ instead of $\sigma$. We can thus conclude that
\[
\frac{\partial}{\partial \tau} \phi_\tau*\phi_\rho*z(x) \leq 0.
\]

If $z$ has more bars we note that the double convolution of a single bar of $z$ extends a distance of $\tau+\rho \leq 2 \rho \leq \omega$ outside of the bar and thus will not influence the value of $\phi_\tau*\phi_\rho*z$ inside other bars of $z$. 
\end{proof}

\begin{lemma}\label{lem:F3firstvar}
Let $z\in \B_{\omega}$, $\phi_\sigma \in \mathcal{K}$, and $u$ a minimiser of $F_3$ over $\B$. Denote by $x_i$ the locations of the interfaces of $u$, with $x_0 < x_1 < \ldots$, then we have for every $i$
\begin{equation}\label{eq:F3firstvar}
\phi_\rho*f_\sigma(x_i) = \frac12 +\phi_\rho*\phi_\rho*\overline u_i(x_i),
\end{equation}
where
\[
\overline u_i := \left\{ \begin{array}{ll} u-\chi_{[x_i, x_{i+1}]} & \text{ if } i \text{ is even, i.e. } x_i \text{ is the left interface of a bar of } u,\\ u-\chi_{[x_{i-1}, x_i]} & \text{ if } i \text{ is odd, i.e. } x_i \text{ is the right interface of a bar of } u.\end{array}\right.
\]
(Note that for $i$ even, $\overline u_i=\overline u_{i+1}$.)

Consequently if $\phi_\sigma = \hat \phi_\sigma$,  
$\sigma \leq \rho \leq \frac{\omega}2$, and $\lambda$, $\rho$ and $\sigma$ satisfy in addition
\begin{equation}\label{eq:sigmarhoomegacondition}
\frac2\lambda + \frac1{15 \rho^2} \Bigl(-\sigma^3 + 5 \rho \sigma^2 + 17 \rho^3\Bigr) < \omega,
\end{equation}
then for every $i$, $x_i$ is the location of an interface of $z$.
\end{lemma}

\begin{figure}[ht]
\centerline{\resizebox{3cm}{!}{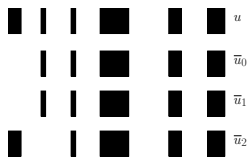}}
\caption{For a given bar code $u$ the bar codes $\overline u_0$, $\overline u_1$, and $\overline u_2$ as used in Lemma~\ref{lem:F3firstvar} are shown}
\label{fig:uubar}
\end{figure}

An example of a bar code $u$ and its accompanying bar codes $\overline u_0$, $\overline u_1$, and $\overline u_2$ is shown in Figure~\ref{fig:uubar}.

\begin{proof}[Proof of Lemma~\ref{lem:F3firstvar}]
Let $\phi_\sigma\in \mathcal{K}$. Let $u$ minimise $F_3$ over $\B$, then $F_3$ has vanishing first variation in $u$ with respect to small perturbations in the locations of the interfaces of $u$. Let $x_0$ be the location of an interface of $u$ where the value of $u$ jumps from $0$ to $1$, in other words, it is the left interface of a bar. The argument is analogous for a right interface. We consider a perturbed $u(t):= \overline u_0 + \chi_{[x_0+t, x_1]} \in \B$, where $|t|$ is small enough such that no interfaces are created or annihilated. The number of interfaces of $u(t)$ is equal to that of $u$ and hence we compute (integration is with respect to $x$)
\begin{align*}
&\hspace{0.7cm}\lambda^{-1} \Big(F_3(u(t)) - F_3(u)\Big) \\
&= \int_\R \bigg[ \Big(\phi_\rho*u(t)-f_\sigma\Big)^2 - \Big(\phi_\rho*u-f_\sigma\Big)^2\bigg]\\
&= \int_\R \bigg[ \Big(\phi_\rho*u(t)\Big)^2 - \Big(\phi_\rho*u\Big)^2 + 2 f_\sigma \cdot \phi_\rho*\Big(u-u(t)\Big)\bigg]\\
&= \left\{ \begin{array}{ll}
\mathlarger\int_\R \bigg[
\Big(\phi_\rho*\chi_{[x_0, x_0+t]}\Big)^2 - 2 \phi_\rho*u \cdot \phi_\rho*\chi_{[x_0, x_0+t]}
+ 2 \phi_\rho*f_\sigma \cdot\Big(u-u(t)\Big)\bigg] & \text{if } t>0,\\
\mathlarger\int_\R \bigg[
\Big(\phi_\rho*\chi_{[x_0+t, x_0]}\Big)^2 + 2 \phi_\rho*u \cdot \phi_\rho*\chi_{[x_0+t, x_0]}
+ 2 \phi_\rho*f_\sigma \cdot\Big(u-u(t)\Big)\bigg] & \text{if } t<0,
\end{array}\right.
\end{align*}
where we have used that $u(t)=u-\chi_{[x_0, x_0+t]}$ if $t>0$ and $u(t)=u+\chi_{[x_0+t, x_0]}$ if $t<0$ in the last line as well as using Lemma~\ref{eq:intconv}.

Assume for now that $t>0$. The case for $t<0$ is analogous. Then, using $u(t)=\overline u_0+\chi_{[x_0+t, x_1]}$,
\begin{align}
\left.\frac{d}{dt} \lambda^{-1} \Big(F_3(u(t)) - F_3(u)\Big)\right|_{t=0^+} &= \left[2 \mathlarger\int_\R \left\{\left(\frac{d}{dt} \int_{x_0}^{x_0+t} \phi_{\rho}(x-y)\,dy\right) \cdot \int_{x_0}^{x_0+t} \phi_\rho(x-y)\, dy\right.\right.\notag\\&\hspace{0.7cm}\left.\left. - 2 \frac{d}{dt} \int_{x_0}^{x_0+t} \phi_{\rho}*\phi_{\rho}*u + 2 \frac{d}{dt} \int_{x_0}^{x_0+t} \phi_\rho*f_\sigma\right\}\,dx \right]_{t=0}\notag\\
&= -2 \phi_\rho*\phi_\rho*u(x_0) - 2 \phi_\rho*f_\sigma(x_0).\label{eq:derivativewrtt}
\end{align}
We can rewrite the first terms as follows:
\[
\phi_\rho*\phi_\rho*u(x_0) = \phi_\rho*\phi_\rho*(\chi_{[x_0, x_1]} + \overline u_0)(x_0) = \frac12 + \phi_\rho*\phi_\rho*\overline u_0(x_0),
\]
where we have used (\ref{eq:equaltoahalf}) to compute $\displaystyle \phi_\rho*\phi_\rho*\chi_{[x_0, x_1]}(x_0) = \frac12$.

Vanishing of the first variation tells us that the right hand side in (\ref{eq:derivativewrtt}) is zero and hence
\[
1 + 2  \phi_\rho*\phi_\rho*\overline u_0(x_0) - 2 \phi_\rho*f_\sigma(x_0) = 0,
\]
which gives equation (\ref{eq:F3firstvar}) for $x_i=x_0$.



Now assume $\sigma \leq \rho \leq \frac{\omega}2$ and $\phi_\sigma=\hat\phi_\sigma$. If $u$ is such that the white spaces between every two subsequent black bars have widths of at least $2\rho$
then it follows that $\hat\phi_\rho*\hat\phi_\rho*\overline u_0(x_i)=0$ for every $i$ and equation (\ref{eq:F3firstvar}) reduces to
\[
\hat\phi_\rho*f_\sigma(x_i) = \frac12.
\]
Lemma~\ref{lem:levelhalfrhosigmaz} then completes the argument. Note that in this case condition (\ref{eq:sigmarhoomegacondition}) is not necessary.

Now assume that $u$ is not as above, i.e. there exist two bars in $u$ separated by a white space of width strictly less than $2\rho$. We will show a contradiction. Let $x_1$ be the right interface of a bar of $u$ and let $x_2$ be the left interface of the next bar, such that $x_2-x_1 < 2\rho \leq \omega$. Then the following inequalities should be satisfied
\[
\hat\phi_\rho*f_\sigma(x_i) = \frac12 + \hat\phi_\rho*\hat\phi_\rho*\overline u_i(x_i) \geq \frac12, \quad \text{for } i\in\{1, 2\}.
\]
According to Lemma~\ref{lem:levelhalfrhosigmaz} this means that $x_1, x_2 \in \supp z$. Now there are two possibilities. The first is that $x_1$ and $x_2$ are located in different bars of $z$. Since $z\in \B_{\omega}$ this means that $x_2-x_1 \geq \omega$ which contradicts our assumption. The second possibility is that $x_1$ and $x_2$ are in the same bar of $z$. Assume the latter now.

By the same arguments the right interface of the second bar, i.e. $x_3$ also lies in $\supp z$. It can lie either in a different bar of $z$ than $x_2$ or in the same one. In the former case we have that there exists an interval $N$ with $|N| \geq \omega$ such that $z=0$ and $u=1$ on $N$ and using (\ref{eq:sigmarhoomegacondition}) we can use the arguments as in Lemma~\ref{lem:F3partresult} to arrive at a contradiction with the fact that $u$ is a minimiser of $F_3$.\footnote{We don't need $u\in\B^{ij}$ here, because we know that in this construction $N$ has a distance of at least $\omega$ to $x=0$ and to $x=1$.} We conclude that $x_2$ and $x_3$ must lie in the same bar of $z$. In a similar way we find that $x_0$ lies in the same bar.
If $z$ has more than two bars, via induction on the interfaces we find that
for every even $i$, $[x_i, x_{i+1}] \subset \supp z$. In words, every bar of $u$ is contained in a bar of $z$.

From the foregoing we deduce that $(u-z)(x) \in \{-1, 0\}$ a.e. and
\begin{equation}\label{eq:uminzestimate}
u-z \leq -\chi_{[x_1, x_2]}.
\end{equation}

Define $\hat u := u + \chi_{[x_1, x_2]}$, then
\begin{align}
&\hspace{0.6cm} \int_\R \bigg( \Big(\hat\phi_\rho*\hat u-\hat\phi_\sigma*z\Big)^2 - \Big(\hat\phi_\rho*u-\hat\phi_\sigma*z\Big)^2 \bigg)\notag\\
&= \int_\R \bigg( \Big(\hat\phi_\rho*\hat u\Big)^2 + 2 \hat\phi_\sigma*z \cdot \hat\phi_\rho*(u-\hat u) - \Big(\hat\phi_\rho*u\Big)^2 \bigg)\notag\\
&= \int_\R \Big(\hat\phi_\rho * \chi_{[x_1, x_2]}\Big)^2 + 2 \int_\R \Big(\hat\phi_\rho*u\cdot \hat\phi_\rho*\chi_{[x_1, x_2]} - \hat\phi_\sigma*z \cdot \hat\phi_\rho*\chi_{[x_1, x_2]}\Big)\notag\\
&= \int_\R \Big(\hat\phi_\rho * \chi_{[x_1, x_2]}\Big)^2 + 2 \int_{x_1}^{x_2} \Big(\hat\phi_\rho*\hat\phi_\rho*u - \hat\phi_\rho*\hat\phi_\sigma*z\Big),\label{eq:F3L2uhatucompare}
\end{align}
where the last equality follows by Lemma~\ref{eq:intconv}.

We now use Lemma~\ref{lem:decreasing} and inequality (\ref{eq:uminzestimate}) to estimate
\begin{align*}
\int_{x_1}^{x_2} \Big(\hat\phi_\rho*\hat\phi_\rho*u - \hat\phi_\rho*\hat\phi_\sigma*z \Big) &\leq \int_{x_1}^{x_2} \hat\phi_\rho*\hat\phi_\rho*(u-z) 
\leq -\int_{x_1}^{x_2} \hat\phi_\rho*\hat\phi_\rho*\chi_{[x_1, x_2]}\\
&=-\int_{x_1}^{x_2} \int_\R \int_{x_1}^{x_2} \hat\phi_\rho(x-y) \hat\phi_\rho(y-q) \,dq \,dy \,dx\\
&= -\int_\R \bigg(\int_{x_1}^{x_2} \hat\phi_\rho(y-x) \,dx\bigg)^2 \,dy 
= -\int_\R \Big(\hat\phi_\rho*\chi_{[x_1, x_2]}\Big)^2.
\end{align*}

Using this in (\ref{eq:F3L2uhatucompare}) we find
\[
F_3(\hat u) - F_3(u) \leq -2 - \lambda \int_\R \Big(\hat\phi_\rho*\chi_{[x_1, x_2]}\Big)^2 < 0
\]
which contradicts $u$ being a minimiser.
\end{proof}

\begin{remark}
The result of Lemma~\ref{lem:F3firstvar} doesn't change if $z\in B_\omega^{ij}$ and we minimise $F_3$ over $B^{ij}$ for $i, j \in \{0, 1\}$. Also note that the result can be obtained for any $\phi_\sigma \in \mathcal{K}_3$ if we replace condition (\ref{eq:sigmarhoomegacondition}) by the corresponding parameter range for that choice of kernel, which we can obtain be redoing the calculations in the proof of Lemma~\ref{lem:F3partresult} after (\ref{eq:F4inequality}) for the new kernel.
\end{remark}

\bigskip

Note that in the case where $\sigma \leq \rho$ we could have used Lemma~\ref{lem:F3firstvar} in the proof of Lemma~\ref{lem:F3partresult} instead of Lemma~\ref{lemma1}.

\bigskip

\noindent{\bf Proof of Theorem~\ref{thm:F4}, part~\ref{item:F3zminoverB}}:
From Lemma~\ref{lem:F3firstvar} it follows that under the stated conditions the only possible minimisers of $F_3$ over $\B^{ij}$ are $u=z$ or a $u$ with strictly less interfaces than $z$.
By Lemma~\ref{lem:F3partresult} however such a minimiser cannot have less interfaces than $z$ and hence by Lemma~\ref{lemma0} $u=z$ is the unique minimiser of $F_3$ over $\B^{ij}$.

\qed

\section{Numerical simulations}\label{numerics} 
We present a few test simulations for the minimisation problems  $F_2$ and $F_3$. 
To this end, there exists an  increasing number of state of the art techniques  concerning  TV-based minimisation. However here we are not attempting to write the most efficient algorithm, we only aim to test whether the parameter regimes we found theoretically are close to optimal or not. Hence we take  the naive approach 
 of using a phase field to approximate the total variation:  That is,   choosing $\epsilon$ small, we 
replace the total variation with  
\[\mathlarger\int_0^1  \left(\epsilon \, |u'|^2 \,\, + \,\, \frac{  u^2(1 -u)^2 }{2\epsilon}\,\, dx \right),\]
and consider the $L^2$ gradient descent of the resulting functional. 
While this technique brings in diffuse interfaces\footnote{For actual implementation, one would 
need to threshold the output of the minimisation process in order to generate a bar code} 
 (i.e. minimisers will no longer be bar codes), it is  well-justified for small $\epsilon$ (c.f. \cite{Braides02}) 
 in that minimisers will be close to minimisers of the original sharp interface problem. 
One problem with this method in higher dimensions is that one tends to get stuck in metastable states, and 
hence this method would not work well for 2D bar codes. 
However, our 1D problem is sufficiently rigid so that the method works well and fairly quickly.  
It takes seconds to run our Python code, and while more direct state of the art methods  would 
be substantially  faster (as would be needed in a practical application), our limited goals are well 
served by the phase field approach.   

For $F_3$ the $L^2$ gradient flow gives the equation 
\begin{equation}\label{pde}
 u_t \, = \, 2\epsilon u_{xx} \, - \,   \frac{1}{\epsilon} W'(u) \, - \,  2 \lambda \, \phi_\rho \, * \, (\phi_\rho * u - f_\sigma ),
 \end{equation}
where\footnote{This choice of constant prefactor in $W$ does not lead to unit surface tension in the sharp interface limit, hence $\lambda$ in the simulations differs by an $\mathcal{O}(1)$ multiplicative factor from the $\lambda$ in the analytical results in this paper.} $W(u) =  \, \frac{u^2(1 -u)^2}{2}$. 
In all of our experiments, a bar code is generated with $X$-dimension $\omega \approx 0.0133$.
Except for Figure~\ref{plot3} (bottom right), 
convolution with the hat function $\hat \phi_\sigma$ is followed by the addition of noise with amplitude $a = 0.1$\footnote{
The added noise was determined as follows: The 400 grid points making up each interval of length $\omega$ were divided into 16 equal groups, each of which was assigned a random number between $-a$ and $a$.}.   
We used $\e=0.0004$ and initial data was always taken to be either $u\equiv 0$ or $u \equiv 1/2$.

\begin{figure}
\centerline{	{\includegraphics[height=2.8in]{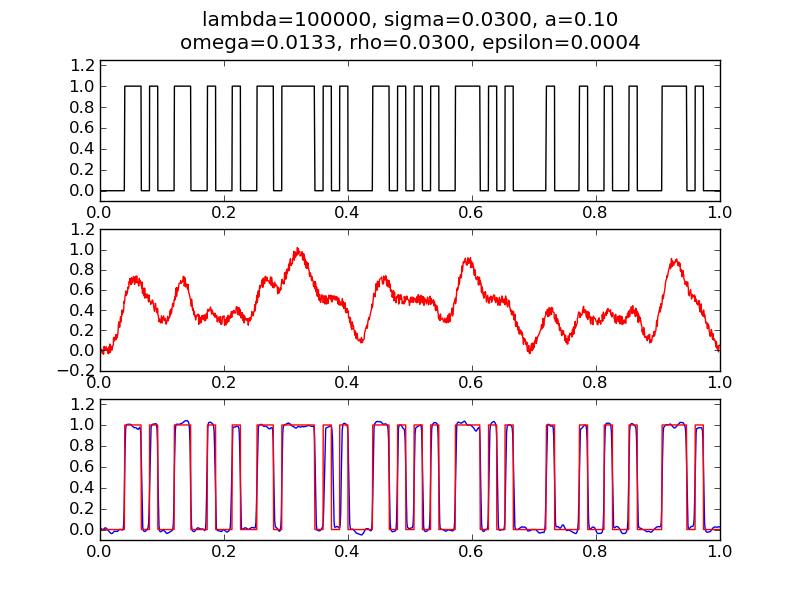}}  \hspace{-1cm} 	{\includegraphics[height=2.8in]{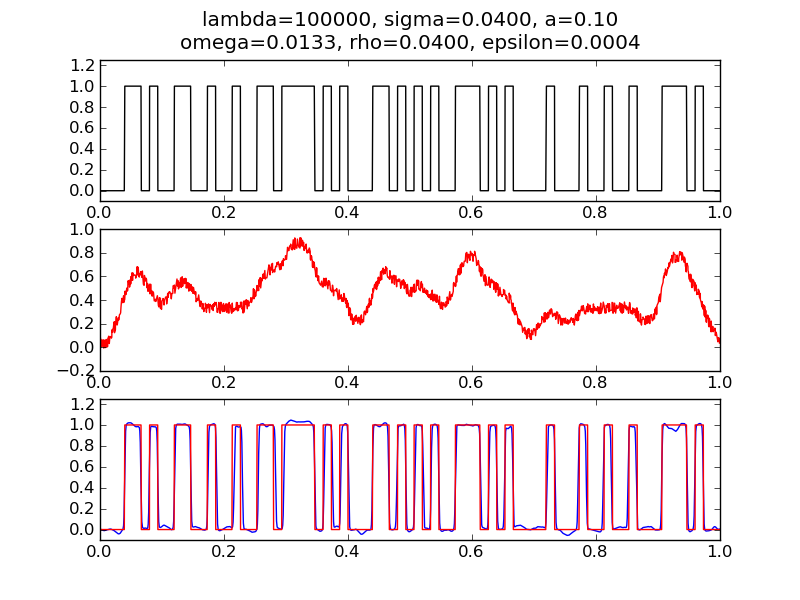}}}
\caption{ Here we look at minimisers of $F_2$ ($\rho = \sigma$) to find that 
the algorithm works for blurring far past $\omega$. 
In all simulations, the three rows are as follows: 
A bar code is generated with $X$-dimension exactly $\omega = 0.0133$;  convolution $f_\sigma$ of the bar code with $\hat \phi_\sigma$ with added noise of amplitude $a$; final steady state for (\ref{pde}) superimposed with the generating bar code.  }
\label{plot1}
\end{figure}

\begin{figure}
\centerline{	{\includegraphics[height=2.8in]{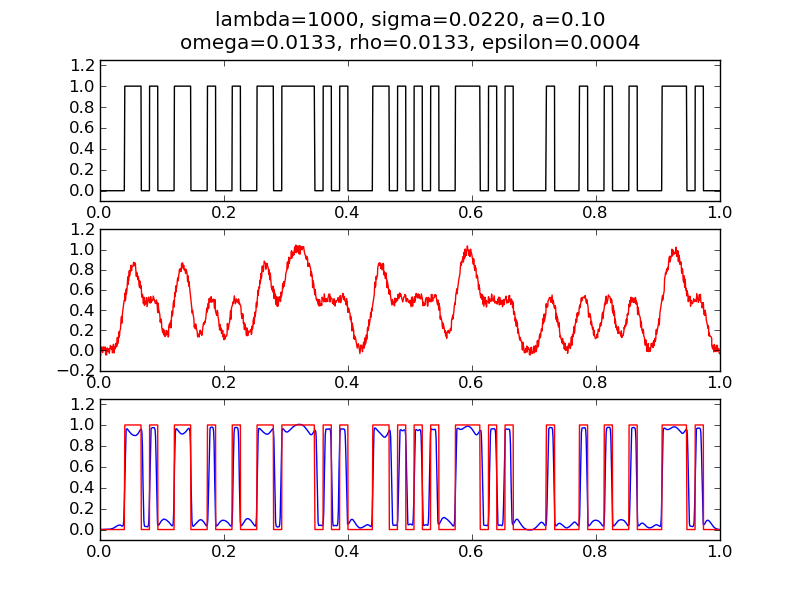}} 	\hspace{-1cm} {\includegraphics[height=2.8in]{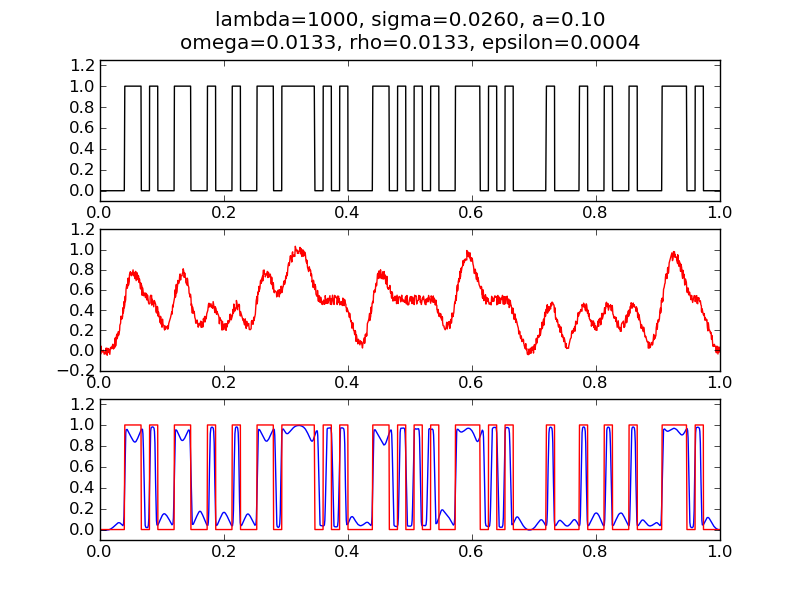}}}
\centerline{	{\includegraphics[height=2.8in]{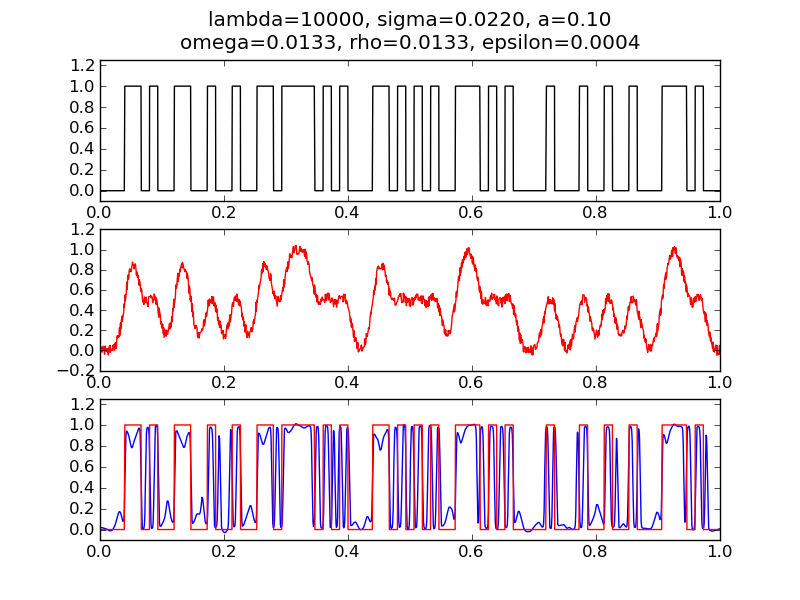}}  \hspace{-1cm} 	{\includegraphics[height=2.8in]{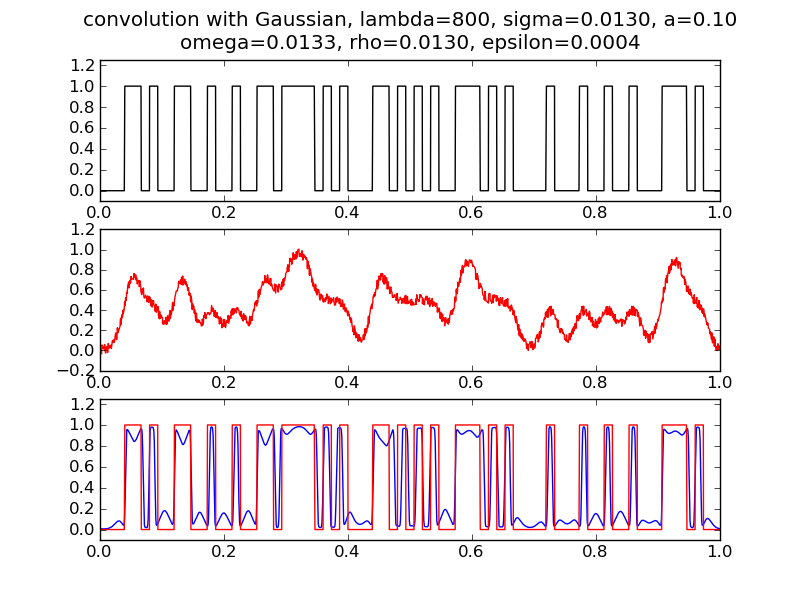}}}
\caption{Minimisers for $F_3$. Here we take the deconvolution kernel size to be the $X$-dimension. Top row: we convolute the data with $X$-dimension exactly $\omega$ for two choices of $\sigma$. The algorithm works well with $\lambda = 1000$. However, as noted in the bottom left, for larger values of $\lambda$ it loses information. Bottom right: a bar code is convoluted with a Gaussian with standard deviation $\sigma$ but deconvoluted with $\hat\phi_\sigma$. }
\label{plot3}
\end{figure}

The algorithm works well for $\sigma$ far beyond the regime of Theorem \ref{thm:F4}. 
We give a few sample results. In 
Figure~\ref{plot1} (left) we see that  choosing  $\rho = \sigma$ (i.e. using $F_2$),  one obtains good results for  $\sigma$ larger than twice $\omega$.   Figure~\ref{plot1} (right) shows that even for $\sigma \approx 3 \omega$, the results are not bad, however they begin to loose accuracy. 
In Figure~\ref{plot3} we note that choosing $\rho$ to be the $X$-dimension works well for blurring  with $\sigma$ up to twice $\omega$. Note that here we are in the regime $\rho  <   \sigma$, which we avoided in Theorem~\ref{thm:F4}. In fact, our 
counter example suggested that in this regime an upper bound on $\lambda$ is necessary. Figure~\ref{plot3} (bottom left)  indeed supports this observation by taking $\lambda$ much larger than in Figure~\ref{plot3} (top right).


We also performed tests where we {\it convolute/blur} the bar code with a Gaussian kernel with  standard deviation $\sigma$ but {\it deconvolute/deblur} with the hat function $\hat \phi_\rho$. To obtain satisfactory results,  one must choose $\sigma$ and $\rho$ very close to each other and no larger than $\omega$, and  use a suitably tuned  midrange  $\lambda$. We give one example in Figure~\ref{plot3} (bottom right). 

Simulations were also performed for $F_1$ (no deconvolution/deblurring  kernel) but we always found that using $F_2$ or $F_3$ (with even a small deblurring kernel) was preferable.  


\section{Discussion}\label{disc}

We have presented results on the accuracy of  TV-based energy minimisation methods for bar code deblurring in certain parameter regimes. Numerical simulations, which included the effects of noise, show that these methods are valid in much larger regimes and in particular, allow for significantly more blurring. While our analytical results did not showcase the benefits of 
using a deconvolution/deblurring kernel (i.e. the merits of  $F_2, F_3$ versus $F_1$), numerical experiments showed clearly that the presence of a  deblurring  kernel  in $F_2$ or  $F_3$ always gave better results over no deblurring ($F_1$). 

In practice,  the size of the blurring kernel pertains to the so-called {\it spot diameter} 
of the laser beam at  impact with the bar code.  This {\it spot diameter}  is a function of the laser beam and the distance from the scanner to the bar code. According to Palmer \cite{Palmer07} (p.127), most scanners can successfully read a bar code if the spot diameter is no greater than $\sqrt{2}$ times the $X-$dimension (i.e. for $2 \sigma < \sqrt{2} \omega$). 
This suggests that the even the conditions we have imposed on $\sigma$ in our results are not completely unreasonable. However, as suggested by the numerics, one might be able to prove results for $\sigma$ past the $X-$dimension. 

Realistically neither the spot diameter (size of the blurring kernel) nor the distribution  of the beam intensity (shape of the blurring kernel) is exactly known,  and inferring this information from signals is an ill-posed problem. 
 In  \cite{Esedoglu04}, the author considers a Gaussian \emph{ansatz} for all kernels but introduces a novel optimization scheme for determining the standard deviation of the blurring kernel. 
In terms of the shape of the kernel, 
our last simulation in Figure~\ref{plot3} (bottom right) is suggestive. We note that if the convolution in the measured signal is done with an infinitely supported Gaussian with standard deviation $\sigma$, then deconvolution with a hat function of approximate size $\sigma$ works reasonably well. Thus if one could determine certain statistics of the blurring kernel, one could then deconvolute with a set kernel possessing similar statistics. 
Determining such statistics should in principle be possible as some standard bar code symbologies have a fixed structure at their left and right boundaries (the {\it left and right guards}, c.f. \cite{Palmer07}).

\bigskip

{\bf Acknowledgments:} This work was completed while both authors were at Simon Fraser University. We are grateful to Fadil Santosa for bringing this problem to our attention and for many interesting conversations. 
We also  thank Selim Esedo\=glu for  useful discussions and for the use of his original code which was the basis for our numerical experiments. This code was modified and tested in \emph{Python} with the \emph{NumPy} package  by Simon Fraser undergraduate student Jacob Groundwater, who we would also like to thank. 
This research  was partially supported by  an 
 NSERC (Canada) Discovery Grant. YvG was also supported by a PIMS postdoctoral fellowship.

\bibliographystyle{acm}
\bibliography{bibliography}

\appendix

\section{Calculations in the proof of Lemma~\ref{lem:F3partresult}}\label{sec:moredetails}

In this appendix we collect some of the longer calculations in the proof of Lemma~\ref{lem:F3partresult}.
We start with a lemma.

\begin{lemma}\label{lem:fsigmasquared}
Let $z:= \chi_{[a, b]}$ for some $a < b$, $\sigma \leq \frac{b-a}2$, and $f_\sigma = \hat\phi_\sigma*z$, then
\[
\int_\R f_\sigma^2 = b-a-\frac7{15} \sigma.
\]
\end{lemma}
\begin{proof}
We compute
\begin{align}\label{eq:fsigmaexplicit1}
f_\sigma(x) &=
\left\{ \begin{array}{ll}
0 & \text{if } x \in (-\infty, a-\sigma],\\
I_+^\sigma(x, a, x+\sigma) & \text{if } x\in [a-\sigma, a],\\
I_-^\sigma(x, a, x) + I_+^\sigma(x, x, x+\sigma) & \text{if } x\in [a, a+\sigma],\\
I_-^\sigma(x, x-\sigma, x) + I_+^\sigma(x, x, x+\sigma) & \text{if } x\in[a+\sigma, b-\sigma],\\
I_-^\sigma(x, x-\sigma, x) + I_+^\sigma(x, x, b) & \text{if } x\in[b-\sigma, b],\\
I_-^\sigma(x, x-\sigma, b) &\text{if } x\in[b, b+\sigma],\\
0 & \text{if } x\in[b+\sigma, \infty).
\end{array}\right. \notag\\
&= \left\{ \begin{array}{ll}
0 & \text{if } x \in (-\infty, a-\sigma],\\
\frac1{2\sigma^2} (x+\sigma-a)^2 & \text{if } x\in [a-\sigma, a],\\
-\frac1{2\sigma^2} \Bigl(x-a-\bigl(1+\sqrt{2}\bigr) \sigma\Bigr) \Bigl(x-a-\bigl(1-\sqrt{2}\bigr) \sigma\Bigr) & \text{if } x\in [a, a+\sigma],\\
1 & \text{if } x\in[a+\sigma, b-\sigma],\\
-\frac1{2\sigma^2} \Bigl(x-b+\bigl(1-\sqrt{2}\bigr)\sigma\Bigr) \Bigl(x-b+\bigl(1+\sqrt{2}\bigr)\sigma\Bigr) & \text{if } x\in[b-\sigma, b],\\
\frac1{2\sigma^2} (x-\sigma-b)^2 &\text{if } x\in[b, b+\sigma],\\
0 & \text{if } x\in[b+\sigma, \infty).
\end{array}\right.
\end{align}
An explicit computation of the integral we are interested in leads to the result.
\end{proof}

\bigskip

Next we give the calculations for the different cases described in the proof of Lemma~\ref{lem:F3partresult}.

\bigskip

\textbf{Case Ia:}
\begin{align}
&\hspace{0.4cm}-2 \int_\R \int_N \hat\phi_\rho(x-y)\,dy \int_{N^c\cap[0,1]} \hat\phi_\sigma(x-w)\,dw\,dx\notag\\
&=-2\Biggl\{ \int_{a-\rho}^{a-\sigma} I_+^\rho(x, a, x+\rho) \biggl(I_-^\sigma(x, x-\sigma, x) + I_+^\sigma(x, x, x+\sigma) \biggr)\,dx\notag\\
&\hspace{1.6cm} +\int_{a-\sigma}^a I_+^\rho(x, a, x+\rho) \biggl(I_-^\sigma(x, x-\sigma, x) + I_+^\sigma(x, x, a)\biggr)\,dx\notag\\
&\hspace{1.6cm} +\int_a^{a+\sigma} \biggl(I_-^\rho(x, a, x) + I_+^\rho(x, x, x+\rho)\biggr) I_-^\sigma(x, x-\sigma, a) \,dx\notag\\
&\hspace{1.6cm} +\int_{a+|N|-\sigma}^{a+|N|} \biggl(I_-^\rho(x, x-\rho, x) + I_+^\rho(x, x, a+|N|) \biggr) I_+^\sigma(x, a+|N|, x+\sigma) \,dx\notag\\
&\hspace{1.6cm} +\int_{a+|N|}^{a+|N|+\sigma} I_-(x, x-\rho, a+|N|) \biggl(I_-^\sigma(x, a+|N|, x) + I_+^\sigma(x, x, x+\sigma) \biggr)\,dx\notag\\
&\hspace{1.6cm} +\int_{a+|N|+\sigma}^{a+|N|+\rho} I_-^\rho(x, x-\rho, a+|N|) \biggl(I_-^\sigma(x, x-\sigma, x) + I_+^\sigma(x, x, x+\sigma) \biggr)\,dx\Biggr\}\notag\\
&= \frac1{15 \rho^2} \Bigl(\sigma^3 - 5 \rho \sigma^2 - 10 \rho^3\Bigr).\label{eq:lemCaseIa}
\end{align}

\bigskip

The way to find the specific intervals of integration in the integrals above (and those that follow below) is similar in spirit to what was done in the proof of Theorem~\ref{thm:minimofF2}. We do not give all the details here, but it is important to reflect on the role of the conditions $z\in \B_\omega^{ij}$, $u_0 \in \B^{ij}$ (instead of $z\in \B_\omega$ and $u\in\B$) and $\rho, \sigma\leq \frac\omega2$. Such considerations are addressed in Remark~\ref{rem:roleofconditions}.

\bigskip

\textbf{Case Ib:}
\begin{align}
&\hspace{0.4cm}-2 \int_\R \int_N \hat\phi_\rho(x-y)\,dy \int_{N^c\cap[0,1]} \hat\phi_\sigma(x-w)\,dw\,dx\notag\\
&=-2 \Biggl\{ \int_{a-\rho}^a I_+^\rho(x, a, x+\rho) \biggl(I_-^\sigma(x, x-\sigma, x) + I_+^\sigma(x, x, a) \biggr)\,dx\notag\\
&\hspace{1.6cm} +\int_a^{a+\rho} \biggl(I_-^\rho(x, a, x) + I_+^\rho(x, x, x+\rho) \biggr) I_-^\sigma(x, x-\sigma, a)\,dx\notag\\
&\hspace{1.6cm} +\int_{a+\rho}^{a+\sigma} \biggl(I_-^\rho(x, x-\rho, x)  + I_+^\rho(x, x, x+\rho) \biggr) I_-^\sigma(x, x-\sigma, a) \,dx\notag\\
&\hspace{1.6cm} +\int_{a+|N|-\sigma}^{a+|N|-\rho} \biggl(I_-^\rho(x, x-\rho, x) + I_+^\rho(x, x, x+\rho) \biggr) I_+^\sigma(x, a+|N|, x+\sigma) \,dx\notag\\
&\hspace{1.6cm} +\int_{a+|N|-\rho}^{a+|N|} \biggl(I_-^\rho(x, x-\rho, x) + I_+^\rho(x, x, a+|N|) \biggr) I_+^\sigma(x, a+|N|, x+\sigma) \,dx\notag\\
&\hspace{1.6cm} +\int_{a+|N|}^{a+|N|+\rho} I_-^\rho(x, x-\rho, a+|N|) \biggl(I_-^\sigma(x, a+|N|, x) + I_+^\sigma(x, x, x+\sigma) \biggr)\,dx\Biggr\}\notag\\
&= \frac1{15\sigma^2} \Bigl(\rho^3 - 5 \sigma \rho^2 - 10 \sigma^3\Bigr).\label{eq:lemCaseIb}
\end{align}

\bigskip

\textbf{Case II, first term:}
\begin{align}
&-2\int_\R \int_N \hat\phi_\rho(x-y)\,dy \int_\R \hat\phi_\rho(u-w)\,dw\,dx\notag\\
&=-2 \Biggl\{ \int_{a-\rho}^a I_+^\rho(x, a, x+\rho) \biggl(I_-^\rho(x, x-\rho, x) + I_+^\rho(x, x, x+\rho)\biggr)\,dx \notag\\
& +\int_a^{a+\rho} \biggl(I_-^\rho(x, a, x)  + I_+^\rho(x, x, x+\rho)\,dy\biggr) 
\cdot 
\biggl(I_-^\rho(x, x-\rho, x) + I_+^\rho(x, x, x+\rho)\biggr)\,dx \notag\\
& +\int_{a+\rho}^{a+|N|-\rho} \biggl(I_-^\rho(x, x-\rho, x) + I_+^\rho(x, x, x+\rho)\biggr) 
\cdot 
 \biggl(I_-^\rho(x, x-\rho, x) + I_+^\rho(x, x, x+\rho) \biggr)\,dx \notag\\
& +\int_{a+|N|-\rho}^{a+|N|} \biggl(I_-^\rho(x, x-\rho, x) + I_+^\rho(x, x, a+|N|)\biggr) 
\cdot 
\biggl(I_-^\rho(x, x-\rho, x) + I_+^\rho(x, x, x+\rho)\biggr) \,dx \notag\\
& +\int_{a+|N|}^{a+|N|+\rho} I_-^\rho(x, x-\rho, a+|N|) \biggl(I_-^\rho(x, x-\rho, x) + I_+^\rho(x, x, x+\rho) \biggr)\,dx\Biggr\} \notag\\
&=-2N.\label{eq:lemCaseIIfirstterm}
\end{align}

\textbf{Case II, second term, IIa:}
\begin{align}
& 2 \int_\R \int_N \hat\phi_\rho(x-y) \,dy \int_N \hat\phi_\sigma(x-w)\,dw\,dx\notag\\
&= 2 \Biggl\{ \int_{a-\sigma}^a I_+^\rho(x, a, x+\rho)  I_+^\sigma(x, a, x+\sigma)\,dx\notag\\
& +\int_a^{a+\sigma} \biggl(I_-^\rho(x, a, x) + I_+^\rho(x, x, x+\rho)\biggr)
\cdot  \biggl(I_-^\sigma(x, a, x) + I_+^\sigma(x, x, x+\sigma) \biggr)\,dx\notag\\
&+\int_{a+\sigma}^{a+\rho} \biggl(I_-^\rho(x, a, x) + I_+^\rho(x, x, x+\rho) \biggr) 
\cdot  \biggl(I_-^\sigma(x, x-\sigma, x) + I_+^\sigma(x, x, x+\sigma) \biggr)\,dx\notag\\
& +\int_{a+\rho}^{a+|N|-\rho} \biggl(I_-^\rho(x, x-\rho, x) + I_+^\rho(x, x, x+\rho)\biggr) 
\cdot  \biggl(I_-^\sigma(x, x-\sigma, x) + I_+^\sigma(x, x, x+\sigma) \biggr) \,dx\notag\\
& +\int_{a+|N|-\rho}^{a+|N|-\sigma} \biggl(I_-^\rho(x, x-\rho, x) + I_+^\rho(x, x, a+|N|)\biggr) 
 \cdot  \biggl(I_-^\sigma(x, x-\sigma, x) + I_+^\sigma(x, x, x+\sigma) \biggr)\,dx\notag\\
& +\int_{a+|N|-\sigma}^{a+|N|} \biggl(I_-^\rho(x, x-\rho, x) + I_+^\rho(x, x, a+|N|) \biggr) 
\cdot  \biggl(I_-^\sigma(x, x-\sigma, x) + I_+^\sigma(x, x, a+|N|)\biggr)\,dx\notag\\
& +\int_{a+|N|}^{a+|N|+\sigma} I_-^\rho(x, x-\rho, a+|N|) I_-^\sigma(x, x-\sigma, a+|N|) \,dx\notag\\ \notag\\
&= \,\,\, 2N + \frac1{15 \rho^2} \Bigl(\sigma^3 - 5 \rho \sigma^2 - 10 \rho^3\Bigr).\label{eq:lemCaseIIa}
\end{align}

\section{Proof that $\hat\phi_\sigma$ satisfies condition (\ref{eq:condpsit})}\label{app:psitaufrhononpos}


In this Appendix we prove that the hat function $\hat\phi_\sigma$ satisfies condition (\ref{eq:condpsit}). We do this by proving a more general result first and then showing that this holds for the hat function in particular.

We use the notation as in (\ref {eq:Ksigma}) and introduce

\begin{lemma}\label{lem:B1}
Use the notation as in (\ref {eq:Ksigma}).
If for each $\tau \in (0, \sigma]$
\begin{enumerate}
\item\label{item:B1case1} either $\displaystyle \frac{\partial}{\partial \tau} p(x, \tau)$ is monotonically increasing in $x$ and $\mathcal{J}(\sigma, \tau, 0, c) \leq 0$ for all $c\geq 2\sigma$,
\item\label{item:B1case2} or $\displaystyle \frac{\partial}{\partial \tau} p(x, \tau)$ is monotonically decreasing in $x$ and $\mathcal{J}(\sigma, \tau, \frac{c}2, c) \leq 0$ for all $c\geq 2\sigma$,
\end{enumerate}
then $\mathcal{J}(\sigma, \tau, x, c) \leq 0$ for all $\tau \in (0, \sigma]$, for all $c \geq 2 \sigma$ and all $x\in [0, c]$, i.e condition (\ref{eq:condpsit}) holds.
\end{lemma}
\begin{proof}
Let $\sigma>0$, $c\geq 2\sigma$, and $\tau\in (0, \sigma]$.

Define $f_\sigma:=\phi_\sigma*\chi_{[0,c]}$  and 
\[
\psi_\tau(x):= \left\{ \begin{array}{ll} \frac{\partial}{\partial \tau} p(-x, \tau)& \text{ if } -\tau\leq x \leq 0,\\ \frac{\partial}{\partial \tau} p(x, \tau) & \text{ if } 0 \leq x \leq \tau, \\ 0 & \text{ otherwise}.  \end{array}\right.
\]
This allows us to rewrite
\[
\mathcal{J}(\sigma, \tau, x, c) = \psi_\tau*f_\sigma(x).
\]

We first consider case~\ref{item:B1case1}. Since $\displaystyle \frac{\partial}{\partial \tau} p(\cdot, \tau)$ is monotonically increasing, the function $-\psi_\tau$ is symmetric unimodal. Since the convolution of two symmetric unimodal functions is again a symmetric unimodal function (see \cite{Uhrin84, EatonPerlman91} and references therein) we find
that $-\psi_\tau*f_\sigma$ is a unimodal function with mode at $\frac{c}2$. Hence $\psi_\tau*f_\sigma(0) = \psi_\tau*f_\sigma(c)$ and $\psi_\tau*f_\sigma(x) \leq \psi_\tau*f_\sigma(0) \leq 0$ for all $x\in[0, c]$.

Next we consider case~\ref{item:B1case2}. In this case $\psi_\tau$ is symmetric unimodal and hence $\psi_\tau*f_\sigma$ is unimodal with mode at $\frac{c}2$ and hence for all $x\in [0, c]$ we have $\psi_\tau*f_\sigma(x) \leq \psi_\tau*f_\sigma(c/2) \leq 0$.
\end{proof}

We complete the proof that $\hat \phi_\sigma$ satisfies condition (\ref{eq:condpsit}) by showing that $\hat\phi_\sigma$ satisfies condition~\ref{item:B1case1} in Lemma~\ref{lem:B1}. Let $\sigma$, $\tau$, $c$ and $x$ satisfy the conditions in (\ref{eq:condpsit}). For the hat function we have $\displaystyle p(x, \tau) = \frac1\tau \left(1-\frac{x}\tau\right)$ and hence $\displaystyle \frac{\partial}{\partial \tau} p(x, \tau) = \frac1{\tau^2}\left(-1+\frac{2x}\tau\right)$ is monotonically increasing in $x$. Furthermore by (\ref{eq:fsigmaexplicit1}) ---with $a=0$, $b=c$--- we find that for $y\in[0, \tau]$ we have $f_\sigma(y)=1-f_\sigma(-y)$. We then compute
\begin{align*}
\mathcal{J}(\sigma, \tau, 0, c) &= \psi_\tau*f_\sigma(0) = \frac1{\tau^2}  \int_0^\tau \left(-1+\frac{2y}\tau\right) \big(f_\sigma(-y)+1-f_\sigma(-y)\big)\,dy\\
&= \frac1{\tau^2} \int_0^\tau \left(-1+\frac{2y}\tau\right)\,dy = 0.
\end{align*}

\end{document}